\documentclass[amssymb,amsmath,11pt]{article}
\usepackage{amssymb,amsmath}
\usepackage{graphicx}

\newcommand{\be}{\begin{equation}}
\newcommand{\ee}{\end{equation}}
\newcommand{\bea}{\begin{eqnarray}}
\newcommand{\eea}{\end{eqnarray}}
\newtheorem{theorem}{Theorem}
\newtheorem{proposition}{Proposition}

\newtheorem{corollary}{Corollary}
\newtheorem{example}{Example}
\newtheorem{remark}{Remark}
\newtheorem{lemma}{Lemma}

\def\1#1{^{(#1)}}
\def\la{\langle}
\def\ra{\rangle}

\begin{document}
\title{Duality Relation for the Hilbert Series\\
of Almost Symmetric Numerical Semigroups}
\author{Leonid G. Fel\\
\\Department of Civil Engineering, Technion, Haifa 3200, Israel\\
\vspace{-.2cm}
\\{\sl e-mail: lfel@tx.technion.ac.il}}
\date{\today}
\maketitle
\def\be{\begin{equation}}
\def\ee{\end{equation}}
\def\bea{\begin{eqnarray}}
\def\eea{\end{eqnarray}}
\def\p{\prime}
\vspace{-.2cm}
\begin{abstract}
We derive the duality relation for the Hilbert series $H\left({\bf d}^m;z
\right)$ of almost symmetric numerical semigroup ${\sf S}\left({\bf d}^m\right)$
combining it with its dual $H\left({\bf d}^m;z^{-1}\right)$. On this basis we 
establish the bijection between the multiset of degrees of the syzygy terms and 
the multiset of the gaps $F_j$, generators $d_i$ and their linear combinations. 
We present the relations for the sums of the Betti numbers of even and odd 
indices separately. We apply the duality relation to the simple case of the 
almost symmetric semigroups of maximal embedding dimension, and give the 
necessary and efficient conditions for minimal set ${\bf d}^m$ to generate such 
semigroups.\\ \\
{\bf Keywords:} Almost symmetric semigroups, the Hilbert series, the Betti 
numbers.\\
{\bf 2000 Mathematics Subject Classification:}  Primary -- 20M14, Secondary -- 
11P81.
\end{abstract}
\newpage
\newpage
\section{Introduction}\label{s1}
This article deals mainly with almost symmetric numerical semigroups which were 
introduced in \cite{bafr97} and present a special class of nonsymmetric 
numerical semigroups ${\sf S}\left({\bf d}^m\right)$ in ${\mathbb N}\cup\{0\}$. 
Throughout the article we assume that ${\sf S}\left({\bf d}^m\right)$ is 
finitely generated by a minimal set of positive integers ${\bf d}^m=\left\{d_1,
\ldots,d_m\right\}$ with finite complement in ${\mathbb N}$, $\#\left\{{\mathbb 
N}\setminus {\sf S}\left({\bf d}^m\right)\right\}<\infty$. We study the
generating function $H\left({\bf d}^m;z\right)$ of such semigroup ${\sf S}\left(
{\bf d}^m\right)$,
\begin{eqnarray}
H\left({\bf d}^m;z\right)=\sum_{s\;\in\;{\sf S}\left({\bf d}^m\right)}z^s\;,
\label{b1}
\end{eqnarray}
which is referred to as {\em the Hilbert series} of ${\sf S}\left({\bf d}^m
\right)$.

Recall the main definitions and known facts on numerical semigroups which are 
necessary here. A semigroup ${\sf S}\left({\bf d}^m\right)=\left\{s\in{\mathbb N
}\cup\{0\}\;|\;s=\sum_{i=1}^m x_i d_i,\;x_i\in {\mathbb N}\cup\{0\}\right\}$, is
said to be generated by {\em minimal set} of $m$ natural numbers $d_1<\ldots<d_
m$, $\gcd(d_1,\ldots,d_m)=1$, if neither of its elements is linearly 
representable by the rest of elements. It is classically known that $d_1\geq m$ 
\cite{heku71} where $d_1$ and $m$ are called {\em the multiplicity} and {\em 
the embedding dimension (edim)} of the semigroup, respectively. If equality 
$d_1=m$ holds then the semigroup ${\sf S}\left({\bf d}^m\right)$ is called of 
maximal edim. {\em The conductor} $c\left({\bf d}^m\right)$ of semigroup ${\sf 
S}\left({\bf d}^m\right)$ is defined by $c\left({\bf d}^m\right):=\min\left\{s
\in{\sf S}\left({\bf d}^m\right)\;|\;s+{\mathbb N}\cup\{0\}\subset {\sf S}\left(
{\bf d}^m\right)\right\}$ and related to {\em the Frobenius number} of 
semigroup, $F\left({\bf d}^m\right)=c\left({\bf d}^m\right)-1$. 

Denote by $\Delta\left({\bf d}^m\right)$ the complement of ${\sf S}\left({\bf 
d}^m\right)$ in ${\mathbb N}$, i.e. $\Delta\left({\bf d}^m\right)={\mathbb N}
\setminus{\sf S}\left({\bf d}^m\right)$, and call it the set of gaps. The 
cardinality ($\#$) of $\Delta\left({\bf d}^m\right)$ is called {\em the genus} 
of ${\sf S}\left({\bf d}^m\right)$, $G\left({\bf d}^m\right):=\#\Delta\left(
{\bf d}^m\right)$. For the set $\Delta\left({\bf d}^m\right)$ introduce the 
generating function $\Phi\left({\bf d}^m;z\right)$ which is related to the 
Hilbert series,
\begin{eqnarray}
\Phi\left({\bf d}^m;z\right)=\sum_{s\;\in\;\Delta\left({\bf d}^m\right)}z^s\;,
\;\;\;\;\;\;\;\;\Phi\left({\bf d}^m;z\right)+H\left({\bf d}^m;z\right)=
\frac1{1-z}\;.\nonumber
\end{eqnarray}
Denote by $t\left({\bf d}^m\right)$ {\em the type} of the numerical semigroup
${\sf S}\left({\bf d}^m\right)$ which coincides with cardinality of set ${\sf S}
^{\prime}\left({\bf d}^m\right)$ that is defined \cite{heku71} as follows,
\begin{eqnarray}
{\sf S}^{\prime}\left({\bf d}^m\right)=\left\{F_j\in {\mathbb Z}\;|\;F_j\not\in
{\sf S}\left({\bf d}^m\right)\;\mbox{and}\;F_j+s\in {\sf S}\left({\bf d}^m
\right),\;\mbox{for}\;\forall\;s\in {\sf S}\left({\bf d}^m\right)\setminus 
\{0\},\;j\leq t\left({\bf d}^m\right)\right\},\label{tr6a}
\end{eqnarray}
and $F_j\neq F_k$ if $j\neq k$. Set ${\sf S}^{\prime}\left({\bf d}^m\right)$ is 
not empty since $F\left({\bf d}^m\right)\in {\sf S}^{\prime}\left({\bf d}^m
\right)$ for any minimal generating set $(d_1,\ldots ,d_m)$. 

The semigroup ring ${\sf k}\left[X_1,\ldots,X_m\right]$ over a field ${\sf k}$ 
of characteristic 0 associated with ${\sf S}\left({\bf d}^m\right)$ is a 
polynomial subring graded by $\deg X_i=d_i$, $i=1,\ldots,m$, and generated by 
all monomials $z^{d_i}$. The Hilbert series $H\left({\bf d}^m;z\right)$ of a 
graded subring ${\sf k}\left[z^{d_1},\ldots,z^{d_m}\right]$ is a rational 
function \cite{stan96}
\begin{eqnarray}
H\left({\bf d}^m;z\right)=\frac{Q\left({\bf d}^m;z\right)}{\prod_{j=1}^m
\left(1-z^{d_j}\right)}\;,\label{tr16}
\end{eqnarray}
where $H\left({\bf d}^m;z\right)$ has a pole $z=1$ of order 1. The numerator
$Q\left({\bf d}^m;z\right)$ is a polynomial in $z$,
\begin{eqnarray}
&&Q\left({\bf d}^m;z\right)=1-Q_1\left({\bf d}^m;z\right)+Q_2\left({\bf d}^m;z
\right)-\ldots+(-1)^{m-1}Q_{m-1}\left({\bf d}^m;z\right)\;,\;\;\;\;
\Sigma_m=\sum_{k=1}^md_k\;,\;\;\;\;\;\;\label{bet05}\\
&&Q_{m-1}\left({\bf d}^m;z\right)=Q_{m-1}^{\p}\left({\bf d}^m;z\right)+z^{F
\left({\bf d}^m\right)+\Sigma_m}\;,\;\;\;\;\deg Q_{m-1}^{\p}\left({\bf d}^m;z
\right)<F\left({\bf d}^m\right)+\Sigma_m\;.\label{bet01}\\
&&Q_i\left({\bf d}^m;z\right)=\sum_{j=1}^{\beta_i\left({\bf d}^m\right)}z^{C_{
j,i}}\;,\;\;\;1\leq i\leq m-1\;,\;\;\;\deg Q_i\left({\bf d}^m;z\right)<\deg Q_{
i+1}\left({\bf d}^m;z\right)\;.\label{bet1}
\end{eqnarray}
In formula (\ref{bet1}) the numbers $C_{j,i}$ and $\beta_i\left({\bf d}^m
\right)$ denote {\em the syzygy degrees} and {\em the Betti numbers}, 
respectively. The summands $z^{C_{j,i}}$ in (\ref{bet1}) stand for the syzygies 
of different kinds and $C_{j,i}$ are the degrees of homogeneous basic invariants
for the syzygies of the $i$th kind,
\begin{eqnarray}
&&C_{j,i}\in {\mathbb N}\;,\;\;\;\;C_{j+1,i}\geq C_{j,i}\;,\;\;\;\;C_{\beta_{i+
1},i+1}>C_{\beta_i,i}\;,\;\;\;\;C_{1,i+1}>C_{1,i}\;,\;\;\;\;\mbox{and}
\nonumber\\
&&C_{j,i}\neq C_{r,i+2k-1}\;,\;\;\;1\leq j\leq \beta_i\left({\bf d}^m\right)\;,
\;\;\;1\leq r\leq \beta_{i+2k-1}\left({\bf d}^m\right)\;,\;\;\;
1\leq k\leq \left\lfloor \frac{m-i}{2}\right\rfloor\;.\;\;\;\label{ar10a}
\end{eqnarray}
The last requirement (\ref{ar10a}) means that all necessary cancellations 
(annihilations) of terms $z^{C_{j,i}}$ in (\ref{bet05}) are already performed. 
However the other equalities, $C_{j,i}=C_{r,i+2k}$ and $C_{j,i}=C_{q,i}$, $j\neq
q$, are not forbidden excluding the syzygy degrees of the last $(m-1)$th kind 
\cite{feai07}. The numbers of terms $z^{C_{j,i}}$ in summands are determined by 
$\beta_i\left({\bf d}^m\right)$ which satisfy the equality \cite{stan96}
\begin{eqnarray}
1-\beta_1\left({\bf d}^m\right)+\beta_2\left({\bf d}^m\right)-\ldots +(-1)^{m-1}
\beta_{m-1}\left({\bf d}^m\right)=0\;.\label{bet2}
\end{eqnarray}
The Betti numbers $\beta_i\left({\bf d}^m\right)$ satisfy also an inequality 
(see \cite{feai07}, Theorem 7),
\begin{eqnarray}
1+\beta_1\left({\bf d}^m\right)+\beta_2\left({\bf d}^m\right)+\ldots+\beta_{
m-1}\left({\bf d}^m\right)\leq d_12^{m-1}-2(m-1)\;.\label{bet22}
\end{eqnarray}
Following \cite{feai07} denote by ${\mathbb B}_i\left({\bf d}^m\right)$ the set 
of degrees of the terms $z^{C_{j,i}}$ (up to degeneration, $C_{j,i}=C_{q,i}$, 
$j\neq q$) for syzygies of the $i$th kind which are entering 
$Q_i\left({\bf d}^m;z\right)$ in (\ref{bet1}),
\begin{eqnarray}
{\mathbb B}_i\left({\bf d}^m\right)=\left\{C_{j,i}\in {\mathbb N}\;|\;z^{C_{j,i}
}\in Q_i\left({\bf d}^m;z\right)\;,\;1\leq j\leq \beta_i\left({\bf d}^m\right)
\right\}\;.\label{bet3}
\end{eqnarray}
A containment $(\in )$ in (\ref{bet3}) means that a monomial $z^{C_{j,i}}$ 
enters polynomial $Q_i\left({\bf d}^m;z\right)$ at least once.
 
By (\ref{ar10a}) we conclude that any two sets ${\mathbb B}_i\left({\bf d}^m
\right)$ and ${\mathbb B}_q\left({\bf d}^m\right)$, whose indices differ by odd 
number, $|q-i|=2k-1$, are disjoined, i.e.
\begin{eqnarray}
{\mathbb B}_i\left({\bf d}^m\right)\bigcap {\mathbb B}_{i+2k-1}\left({\bf d}^m
\right)=\emptyset\;,\;\;\;1\leq i\leq m-1\;,\;\;\;
1\leq k\leq \left\lfloor \frac{m-i}{2}\right\rfloor\;.\label{ar10b}
\end{eqnarray}
Note that the relation ${\mathbb B}_i\left({\bf d}^m\right)\bigcap{\mathbb B}_{
i+2k}\left({\bf d}^m\right)\neq\emptyset$ is not forbidden. Let $\oplus$ denote 
{\em a sumset} of the finite set ${\mathbb U}\subset {\mathbb N}$ of integers 
$u_p$ with an integer $\alpha$, ${\mathbb U}\oplus\{\alpha\}=\left\{u_p+\alpha
\;|\;u_p\in{\mathbb U}\right\}$. Then we have,
\begin{lemma}\label{lem1}{\rm (\cite{feai07}, Lemma 1)}
The following equality holds
\begin{eqnarray}
{\mathbb B}_{m-1}\left({\bf d}^m\right)={\sf S}^{\prime}\left({\bf d}^m\right)
\oplus\left\{\Sigma_m\right\}\;.\label{bet3a}
\end{eqnarray}
\end{lemma}
By consequence of (\ref{tr6a}) and (\ref{bet3a}) we get the known equality 
\cite{heku71}, $\beta_{m-1}\left({\bf d}^m\right)=t\left({\bf d}^m\right)$. 
\section{Two Sorts of Gaps in Numerical Semigroups}\label{s2}
Following \cite{jag77} decompose the set of gaps $\Delta\left({\bf d}^m\right)$ 
into two sets $\Delta_{{\cal G}}\left({\bf d}^m\right)$ and $\Delta_{{\cal H}}
\left({\bf d}^m\right)$,
\begin{eqnarray}
&&\Delta_{{\cal G}}\left({\bf d}^m\right)=\left\{g\not\in {\sf S}\left({\bf
d}^m\right)\;|\;F\left({\bf d}^m\right)-g\in {\sf S}\left({\bf d}^m\right)
\right\},\;\;\;\#\Delta_{{\cal G}}\left({\bf d}^m\right)=c\left({\bf
d}^m\right)-G\left({\bf d}^m\right),\;\;\;\;\;\;\;\;\;\;\;\;\nonumber\\
&&\Delta_{{\cal H}}\left({\bf d}^m\right)=\left\{h\not\in {\sf S}\left({\bf
d}^m\right)\;|\;F\left({\bf d}^m\right)-h\not\in {\sf S}\left({\bf d}^m\right)
\right\},\;\;\;\#\Delta_{{\cal H}}\left({\bf d}^m\right)=2G\left({\bf d}^m
\right)-c\left({\bf d}^m\right).\;\;\;\;\;\;\;\;\;\;\;\;\label{tr9}
\end{eqnarray}
The following Theorem is essential in this article.
\begin{theorem}\label{the1}{\rm  (\cite{feai07}, Theorem 1)} 
Let the numerical semigroup ${\sf S}\left({\bf d}^m\right)$ be given with its 
Hilbert series $H\left({\bf d}^m;z\right)$. Then the generating functions for 
the sets $\Delta_{{\cal H}}\left({\bf d}^m \right)$ and $\Delta_{{\cal G}}\left(
{\bf d}^m\right)$ are given by
\begin{eqnarray}
\sum_{h\;\in\;\Delta_{{\cal H}}\left({\bf d}^m\right)}z^h&=&-H\left({\bf d}^m;
z\right)-H\left({\bf d}^m;z^{-1}\right)\cdot z^{F\left({\bf d}^m\right)}\;,
\label{tr17a}\\
\sum_{g\;\in\;\Delta_{{\cal G}}\left({\bf d}^m\right)}z^g&=&\frac{1}{1-z}+
H\left({\bf d}^m;z^{-1}\right)\cdot z^{F\left({\bf d}^m\right)}\;.\nonumber
\end{eqnarray}
\end{theorem}
The next Lemma establishes relationship between the sets, ${\sf S}^{\prime}
\left({\bf d}^m\right)$ and $\Delta_{{\cal H}}\left({\bf d}^m\right)$ for 
general numerical semigroups, and also gives a basis for definition of symmetric
and almost symmetric semigroups.
\begin{lemma}\label{lem2}{\rm (\cite{bafr97} and \cite{feai07}, Lemma 5)} 
Let a numerical semigroup ${\sf S}\left({\bf d}^m\right)$ be given. Then
\begin{eqnarray}
{\sf S}^{\prime}\left({\bf d}^m\right)\setminus\left\{F\left({\bf d}^m\right)
\right\}\subseteq \Delta_{{\cal H}}\left({\bf d}^m\right)\;.\label{tr18}
\end{eqnarray}
\end{lemma}
\subsection{Symmetric, pseudosymmetric and almost symmetric semigroups}
\label{s21}
Imposing requirements on the set $\Delta_{{\cal H}}\left({\bf d}^m\right)$ one 
can simplify significantly the structure of semigroup. A semigroup ${\sf S}
\left({\bf d}^m\right)$ is called {\em symmetric} if $\Delta_{{\cal H}}\left(
{\bf d}^m\right)=\emptyset$ that by (\ref{tr18}) implies $t\left({\bf d}^m
\right)=1$. By Theorem \ref{the1} the following duality relation holds for 
symmetric semigroups (see \cite{feai07}, Corollary 1),
\begin{eqnarray}
H\left({\bf d}^m;z\right)+H\left({\bf d}^m;z^{-1}\right)\cdot z^{F({\bf 
d}^m)}=0\;.\label{tr19a}  
\end{eqnarray}
Notably, all semigroups ${\sf S}\left({\bf d}^2\right)$ are symmetric. For 
$m\geq 3$ the necessary conditions for the minimal set ${\bf d}^m$ to generate 
a symmetric semigroup were given in \cite{wata73}, Lemma 1.

Another simplification comes if $F\left({\bf d}^m\right)$ is an even number, and
$\Delta_{{\cal H}}\left({\bf d}^m\right)=\left\{F\left({\bf d}^m\right)/2
\right\}$. These semigroups are called {\em pseudosymmetric} and the 
corresponding duality relation reads
\begin{eqnarray}
H\left({\bf d}^m;z\right)+H\left({\bf d}^m;z^{-1}\right)\cdot z^{F({\bf 
d}^m)}+z^{\frac{1}{2}F({\bf d}^m)}=0\;.\label{tr19b}
\end{eqnarray}
Pseudosymmetric semigroups have necessarily $t\left({\bf d}^m\right)=2$, but 
the opposite statement (sufficient condition) is not true. For $m=3$, the 
structure of the minimal triple ${\bf d}^3$ generating a pseudosymmetric 
semigroup was given independently in \cite{ros05} and \cite{feai07}, Theorem 9.
Note that both relations, (\ref{tr19a}) and (\ref{tr19b}), are self-dual under 
transformation $z\rightarrow z^{-1}$.

The {\em almost symmetric semigroups} were introduced in \cite{bafr97} as a
generalization of the symmetric and pseudosymmetric ones. They are defined by 
a set equality in (\ref{tr18}), 
\begin{eqnarray}
{\sf S}^{\prime}\left({\bf d}^m\right)\setminus\left\{F\left({\bf d}^m\right)
\right\}=\Delta_{{\cal H}}\left({\bf d}^m\right)\;.\label{tr20}
\end{eqnarray}
Equivalence of (\ref{tr20}) and another equality, $t\left({\bf d}^m\right)=1+
\#\Delta_{{\cal H}}\left({\bf d}^m\right)$, was proven in \cite{bafr97}. For
$m=3$, almost symmetric and pseudosymmetric semigroups coincide.  A minimal set 
${\bf d}^m$ of special kind generating an almost symmetric semigroup is 
given in the next Proposition.
\begin{proposition}\label{pro1}{\rm (\cite{bafr97}, Proposition 11)}
Let a numerical semigroup ${\sf S}\left({\bf d}^{t+1}\right)$ of maximal edim is
generated by tuple $(t+1,t+1+\frac{g}{t},t+1+2\frac{g}{t},\ldots ,t+1+g)$, where
$t\geq 1$, $g\geq -1$ such that $t\mid g$ and $\gcd(t+1,\frac{g}{t})=1$. Then 
${\sf S}\left({\bf d}^{t+1}\right)$ is almost symmetric semigroup with ${\sf S}
^{\prime}\left({\bf d}^{t+1}\right)=\left\{\frac{g}{t},2\frac{g}{t},\ldots ,g
\right\}$.
\end{proposition}

In section \ref{s7} we consider the almost symmetric semigroups with maximal 
edim of generic kind (not satisfying Proposition \ref{pro1}) and give the
necessary and efficient conditions for minimal set ${\bf d}^m$ to generate 
such semigroups.
\section{Duality Relation for Almost Symmetric Semigroups}\label{s3}
Continuing a similar description of symmetric, pseudosymmetric and almost
symmetric semigroups we derive here the duality relation for the Hilbert series 
$H\left({\bf d}^m;z\right)$ for the last ones. By Lemma \ref{lem1} and 
definition (\ref{tr20}) we have 
\begin{eqnarray}
{\mathbb B}_{m-1}\left({\bf d}^m\right)=\left[\Delta_{{\cal H}}\left({\bf d}^m
\right)\cup\left\{F\left({\bf d}^m\right)\right\}\right]\oplus\left\{\Sigma_m
\right\}\;.\label{tr22}
\end{eqnarray}
Following \cite{fel04}, introduce two functions, $\tau$ and its inverse $\tau^{
-1}$, where $\tau$ maps each polynomial $\Psi(z)=\sum c_kz^k\in{\mathbb N}[z]$ 
with $c_k\in\{0,1\}$ onto the set of degrees ${\mathbb K}=\{k\in {\mathbb N}\;
\bracevert\;c_k\neq 0\}$. Since all coefficients of the polynomial $\Psi(z)$ are
1 or 0, we can uniquely reconstruct a set ${\mathbb K}$ and vice versa. In this 
sense $\tau$ is an isomorphic map. The map $\tau$ is also linear in the 
following sense (see \cite{fel04}, Ch. 5):
\begin{eqnarray}
\mbox{If}\;\;\;\;{\mathbb U}_1,{\mathbb U}_2\subset {\mathbb N}\;,\;\;
{\mathbb U}_1\bigcap {\mathbb U}_2=\emptyset\;,\;\;\;\;\mbox{then}\;\;\;\;
\tau^{-1}\left[{\mathbb U}_1\bigcup {\mathbb U}_2\right]=\tau^{-1}\left[
{\mathbb U}_1\right]+\tau^{-1}\left[{\mathbb U}_2\right]\;.\label{tr22a}
\end{eqnarray}
In particular, by consequence of (\ref{bet3}) and Lemma \ref{lem1} we have 
\begin{eqnarray}
\tau^{-1}\left[{\mathbb B}_{m-1}\left({\bf d}^m\right)\right]=Q_{m-1}\left(
{\bf d}^m;z\right)\;.\label{tr22b}
\end{eqnarray}
\begin{theorem}\label{the2}
Let a semigroup ${\sf S}\left({\bf d}^m\right)$ be almost symmetric. Then its 
Hilbert series $H\left({\bf d}^m;z\right)$ satisfies the duality relation
\begin{eqnarray}
H\left({\bf d}^m;z\right)+H\left({\bf d}^m;z^{-1}\right)\cdot z^{F\left({\bf d}
^m\right)}+Q_{m-1}^{\p}\left({\bf d}^m;z\right)\cdot z^{-\Sigma_m}=0\;.
\label{tr24a}
\end{eqnarray}
where $Q_{m-1}^{\p}\left({\bf d}^m;z\right)$ is defined in (\ref{bet01}).
\end{theorem}
{\sf Proof} $\;\;\;$Acting on the left hand side ({\em l.h.s.}) and right 
hand side ({\em r.h.s.}) of (\ref{tr22}) by $\tau^{-1}$ and applying 
(\ref{tr22a}) and (\ref{tr22b}) we get
\begin{eqnarray}
z^{-\Sigma_m}\cdot Q_{m-1}\left({\bf d}^m;z\right)=z^{F\left({\bf d}^m\right)}+
\sum_{h\;\in\;\Delta_{{\cal H}}\left({\bf d}^m\right)}z^h\;.\label{tr23}
\end{eqnarray}
Combining (\ref{tr23}) with (\ref{tr17a}) we obtain,
\begin{eqnarray}
H\left({\bf d}^m;z\right)+H\left({\bf d}^m;z^{-1}\right)\cdot z^{F\left({\bf 
d}^m\right)}=z^{F\left({\bf d}^m\right)}-Q_{m-1}\left({\bf d}^m;z\right)\cdot 
z^{-\Sigma_m}\;.\nonumber
\end{eqnarray}
Substituting the representation (\ref{bet01}) for $Q_{m-1}\left({\bf d}^m;z
\right)$ into the last equation we come to the duality relation (\ref{tr24a}) 
for the Hilbert series for the almost symmetric semigroups.$\;\;\;\;\;\;\Box$

Require that relation (\ref{tr24a}) be self-dual under transformation $z
\rightarrow z^{-1}$, apply it to (\ref{tr24a}) and get 
\begin{eqnarray}
H\left({\bf d}^m;z^{-1}\right)+H\left({\bf d}^m;z\right)\cdot z^{-F\left({\bf 
d}^m\right)}+Q_{m-1}^{\p}\left({\bf d}^m;z^{-1}\right)\cdot z^{\Sigma_m}=0\;.
\label{tr24b}
\end{eqnarray}
By comparison of (\ref{tr24a}) and (\ref{tr24b}) we obtain a necessary condition
for ${\sf S}\left({\bf d}^m\right)$ to be almost symmetric, 
\begin{eqnarray}
Q_{m-1}^{\p}\left({\bf d}^m;z\right)\cdot z^{-\Sigma_m}=Q_{m-1}^{\p}\left({\bf 
d}^m;z^{-1}\right)\cdot z^{F\left({\bf d}^m\right)+\Sigma_m}\;.\label{tr24c}
\end{eqnarray}
The last equation has clear explanation. Indeed, consider $Q_{m-1}^{\p}\left(
{\bf d}^m;z\right)$ which has acording to (\ref{bet01}) and Lemma \ref{lem1} the
following form, $Q_{m-1}^{\p}\left({\bf d}^m;z\right)=\sum_{j=1}^{t\left({\bf 
d}^m\right)-1}z^{F_j+\Sigma_m}$, and substitute it into (\ref{tr24c}),
\begin{eqnarray}
\sum_{j=1}^{t\left({\bf d}^m\right)-1}z^{F_j}=\sum_{j=1}^{t\left({\bf d}^m
\right)-1}z^{-F_j+F\left({\bf d}^m\right)}\;.\label{tr31}
\end{eqnarray}
Formula (\ref{tr31}) is equivalent to the following sequence of equalities,
\begin{eqnarray}
F_{u_1}=F\left({\bf d}^m\right)-F_{v_1}\;,\;\;F_{u_2}=F\left({\bf d}^m\right)-
F_{v_2}\;,\ldots ,\;\;F_{u_{t-1}}=F\left({\bf d}^m\right)-F_{v_{t-1}}\;,
\nonumber
\end{eqnarray}
where $u_1,\ldots ,u_{t-1}$ and $v_1,\ldots ,v_{t-1}$ account for different 
arrangements of the set ${\sf S}^{\prime}\left({\bf d}^m\right)\setminus\left\{
F\left({\bf d}^m\right)\right\}$,
\begin{eqnarray}
\left\{F_{u_1},\ldots ,F_{u_{t-1}}\right\}\equiv \left\{F_{v_1},\ldots ,F_{v_{
t-1}}\right\}\equiv {\sf S}^{\prime}\left({\bf d}^m\right)\setminus\left\{F
\left({\bf d}^m\right)\right\}\;,\;\;\;1\leq u_i,v_i\leq t\left({\bf d}^m\right)
-1\;.\nonumber
\end{eqnarray}
\begin{corollary}\label{cor1}A bijection $\left\{F_{u_1},\ldots,F_{u_{t-1}}
\right\}\Longleftrightarrow\left\{F_{u_1},\ldots ,F_{u_{t-1}}\right\}$ does have
one fixed pont $\frac{1}{2}F\left({\bf d}^m\right)$ iff $t\left({\bf d}^m
\right)$ is an even number, and has not fixed pont iff $t\left({\bf d}^m\right)$
is an odd number.
\end{corollary}
One can make one step further and find the duality relation for the numerator 
$Q\left({\bf d}^m;z\right)$ of the Hilbert series of the almost symmetric 
semigroups. This seams to be reasonable because of the lower syzygies terms 
$z^{C_{j,i}}$, $i\leq m-2$, which are left behind the relation (\ref{tr24c}). 
\begin{theorem}\label{the3}
Let the numerical semigroup ${\sf S}\left({\bf d}^m\right)$ be almost symmetric.
Then the numerator $Q\left({\bf d}^m;z\right)$ of its Hilbert series satisfies 
the duality relation,
\begin{eqnarray}
Q\left({\bf d}^m;z\right)+(-1)^mQ\left({\bf d}^m;z^{-1}\right)\cdot z^{F\left(
{\bf d}^m\right)+\Sigma_m}+Q_{m-1}^{\p}\left({\bf d}^m;z\right)\cdot z^{-\Sigma
_m}\prod_{j=1}^m\left(1-z^{d_j}\right)=0\;.\label{tr34}
\end{eqnarray}
\end{theorem}
{\sf Proof} $\;\;\;$Make use of representation (\ref{tr16}) for the Hilbert 
series $H\left({\bf d}^m;z\right)$ and insert it into (\ref{tr24a}). Multiplying
the r.h.s. and l.h.s. of obtained equation by $\prod_{j=1}^m\left(1-z^{d_j}
\right)$ and keeping in mind the identity $\prod_{j=1}^m\left(1-z^{d_j}\right)/
\left(1-z^{-d_j}\right)=(-1)^mz^{\Sigma_m}$ we arrive at
\begin{eqnarray}
Q\left({\bf d}^m;z\right)+(-1)^mQ\left({\bf d}^m;z^{-1}\right)\cdot z^{F\left(
{\bf d}^m\right)+\Sigma_m}=\prod_{j=1}^m\left(1-z^{d_j}\right)\left[z^{F\left(
{\bf d}^m\right)}-Q_{m-1}\left({\bf d}^m;z\right)\cdot z^{-\Sigma_m}\right]\;.
\label{z32}
\end{eqnarray}
Substituting (\ref{bet01}) for $Q_{m-1}\left({\bf d}^m;z\right)$ into the r.h.s.
of (\ref{z32}) and simplifying it we get (\ref{tr34}).$\;\;\;\;\;\;\Box$

For short, following \cite{jag77} denote $\#\Delta_{{\cal H}}\left({\bf d}^m
\right)=\gamma\left({\bf d}^m\right)$ and by consequence of (\ref{bet05}) write 
another form of (\ref{tr34}) which is more suitable to deal with,
\begin{eqnarray}
\sum_{r=1}^{\gamma\left({\bf d}^m\right)}z^{F_r}\cdot\left(\underbrace{1}-
\sum_{j=1}^mz^{d_j}+\sum_{j>k=1}^mz^{d_j+d_k}-\ldots-(-1)^m\sum_{j=1}^mz^{
\Sigma_m-d_j}+\underbrace{(-1)^mz^{\Sigma_m}}\right)+\;\;\;\;\;\;\nonumber\\
\underbrace{1}-\sum_{j=1}^{\beta_1\left({\bf d}^m\right)}z^{C_{j,1}}+  
\sum_{j=1}^{\beta_2\left({\bf d}^m\right)}z^{C_{j,2}}-\ldots +(-1)^{m-1}\left(
\sum_{r=1}^{\gamma\left({\bf d}^m\right)}\underbrace{z^{F_r+\Sigma_m}}+
\underbrace{z^{F\left({\bf d}^m\right)+\Sigma_m}}\right)+\;\;\;\nonumber\\ 
(-1)^m\left[\underbrace{z^{F\left({\bf d}^m\right)+\Sigma_m}}-\sum_{j=1}^{
\beta_1\left({\bf d}^m\right)}z^{-C_{j,1}+F\left({\bf d}^m\right)+\Sigma_m}+
\ldots+(-1)^{m-1}\left(\sum_{r=1}^{\gamma\left({\bf d}^m\right)}\underbrace{z^{
F\left({\bf d}^m\right)-F_r}}+\underbrace{1}\right)\right]=0\;.\nonumber
\end{eqnarray}
In the last equation we have underbraced $4\gamma\left({\bf d}^m\right)+4$ terms
which are cancelling pairwise. After this simplification and further recasting 
of the rest of terms we get finally,
\begin{eqnarray}
\sum_{r=1}^{\gamma\left({\bf d}^m\right)}\left[\sum_{j=1}^mz^{d_j+
F_r}-\sum_{j>k=1}^mz^{d_j+d_k+F_r}+\ldots+(-1)^m\sum_{j=1}^mz^{\Sigma_m-d_j+F_r}
\right]+\;\;\;\;\;\;\;\;\;\;\;\;\;\label{tr34b}\\
\sum_{j=1}^{\beta_1\left({\bf d}^m\right)}z^{C_{j,1}}-
\sum_{j=1}^{\beta_2\left({\bf d}^m\right)}z^{C_{j,2}}+\sum_{j=1}^{\beta_3\left(
{\bf d}^m\right)}z^{C_{j,3}}-\ldots+(-1)^{m-1}\sum_{j=1}^{\beta_{m-2}\left(
{\bf d}^m\right)}z^{C_{j,m-2}}=\nonumber\\
\;\;\;\;\;(-1)^{m-1}\left[
\sum_{j=1}^{\beta_1\left({\bf d}^m\right)}z^{{\overline C}_{j,1}}-\sum_{j=1}^
{\beta_2\left({\bf d}^m\right)}z^{{\overline C}_{j,2}}+\ldots+(-1)^{m-1}\sum_{
j=1}^{\beta_{m-2}\left({\bf d}^m \right)}z^{{\overline C}_{j,m-2}}\right]\;.
\nonumber
\end{eqnarray}
The degrees ${\overline C}_{j,i}=-C_{j,i}+F\left({\bf d}^m\right)+\Sigma_m$ 
compose another numerical set ${\overline{\mathbb B}}_i\left({\bf d}^m\right)$ 
similar to ${\mathbb B}_i\left({\bf d}^m\right)$,
\begin{eqnarray}
{\overline{\mathbb B}}_i\left({\bf d}^m\right)=\left\{{\overline C}_{j,i}\in 
{\mathbb N}\;|\;z^{C_{j,i}}\in Q_i\left({\bf d}^m;z\right),\;1\leq j\leq\beta_i
\left({\bf d}^m\right)\right\}\;.\label{z35}
\end{eqnarray}
Equation (\ref{tr34b}) is a master equation for study the structure of syzygies 
for the almost symmetric semigroups. In the rest of this paper the main problem 
is considered: what can be said about indeterminate degrees $C_{j,i}$ in Eq.
(\ref{tr34b}), or more precisely, how to find them in terms of the given 
generators $d_i$ and gaps $F_j\in\Delta_{{\cal H}}\left({\bf d}^m\right)$.

A naive approach shows that, in order to satisfy Eq. (\ref{tr34b}), we have
to recast its terms in such a way that in every its l.h.s. and r.h.s. would 
remain only positive terms, and after that to equate all degrees of the terms 
in new recasting equation in its l.h.s. and r.h.s.. In other words, we have to 
build a bijection between two multisets, a multiset ${\mathfrak M}$ of given 
gaps $F_j$, generators $d_i$ and their linear combinations, and a multiset 
${\mathfrak X}$ of indeterminate degrees $C_{j,i}$ and ${\overline C}_{j,i}$.
\section{Multisets and Multiset Operations}\label{s4}
Making preliminary preparation we start with concept of a multiset and basic 
multiset operations allowed by this structure (see \cite{hi80}, \cite{ges95} and
references therein). More rigorous analysis of Eq. (\ref{tr34b}) and associated 
multisets follows later in this and next sections \ref{s5} and \ref{s6}. 

A concept of multiset is a generalization of the concept of a set. A member of 
a multiset can have more than one occurrence (called multiplicity, don't 
confuse with multiplicity of semigroup), while each member of a set occurs once.
A priviledged role is still given to (ordinary) sets when defining maps, 
as there is no clear notion of maps (functions) between multisets.

Consider a multiset which can be formally defined as a pair $\la{\mathbb K},
\sigma_{{\mathbb K}}\ra$ where ${\mathbb K}$ is some finite set, $\#{\mathbb K}<
\infty$, and $\sigma_{{\mathbb K}}:\;{\mathbb K}\mapsto{\mathbb N}$ is a 
function from ${\mathbb K}$ to the set of positive integers. We call $\la
{\mathbb K},\sigma_{{\mathbb K}}\ra$ {\em a standard representation} of 
multiset. For each $\omega\in{\mathbb K}$ the multiplicity (that is, number of 
occurrences) of $\omega$ is the number $\sigma_{{\mathbb K}}(\omega)\geq 0$ such
that
\begin{eqnarray}
\#\la{\mathbb K},\;\sigma_{{\mathbb K}}\ra=\sum_{\omega\in {\mathbb K}}\sigma_{
{\mathbb K}}(\omega)\;,\;\;\;\mbox{and by definition}\;\;\;\mbox{if}\;\;\omega
\not\in{\mathbb K}\;\;\;\mbox{then}\;\;\;\sigma_{{\mathbb K}}(\omega)=0\;.
\label{z8}
\end{eqnarray}
We say that element $\left(\omega,\sigma_{{\mathbb K}}(\omega)\right)$ belongs
to multiset $\la{\mathbb K},\sigma_{{\mathbb K}}\ra$ iff $\sigma_{{\mathbb K}}(
\omega)$ takes a positive value, i.e.
\begin{eqnarray}
\left(\omega,\sigma_{{\mathbb K}}(\omega)\right)\in\la{\mathbb K},\sigma_{
{\mathbb K}}\ra\;\;\mbox{iff}\;\;\sigma_{{\mathbb K}}(\omega)>0\;,\;\;\;\mbox{
and}\;\;\;\left(\omega,\sigma_{{\mathbb K}}(\omega)\right)\not\in\la{\mathbb K}
,\sigma_{{\mathbb K}}\ra\;\;\mbox{iff}\;\;\sigma_{{\mathbb K}}(\omega)=0\;.
\label{z6}
\end{eqnarray}

Let two multisets $\la{\mathbb K}_1,\sigma_{{\mathbb K}_1}\ra$ and $\la{\mathbb 
K}_2,\sigma_{{\mathbb K}_2}\ra$ be given with functions $\sigma_{{\mathbb K}_i}
:\;{\mathbb K}_i\mapsto{\mathbb N}$, $i=1,2$. We say also that the following 
multiset containment ( $\sqsubseteq$ ) holds
\begin{eqnarray}
\la{\mathbb K}_1,\sigma_{{\mathbb K}_1}\ra\sqsubseteq\la{\mathbb K}_2,\sigma_{
{\mathbb K}_2}\ra\;,\;\;\;\mbox{if}\;\;\;{\mathbb K}_1\subseteq{\mathbb K}_2\;
\;\;\mbox{and}\;\;\;\sigma_{{\mathbb K}_1}(\omega)\leq\sigma_{{\mathbb K}_2}(
\omega)\;,\;\;\;\mbox{for}\;\;\forall\;\;\omega\in {\mathbb K}_1,{\mathbb K}_2
\;.\label{bii1}
\end{eqnarray}
In particular, define the multiset equality,
\begin{eqnarray}
\la{\mathbb K}_1,\sigma_{{\mathbb K}_1}\ra=\la{\mathbb K}_2,\sigma_{{\mathbb K}
_2}\ra\;,\;\;\;\mbox{if}\;\;\;{\mathbb K}_1={\mathbb K}_2\;\;\;\mbox{and}\;\;\;
\sigma_{{\mathbb K}_1}(\omega)=\sigma_{{\mathbb K}_2}(\omega)\;,\;\;\;\mbox{for}
\;\;\forall\;\;\omega\in {\mathbb K}_1,{\mathbb K}_2\;.\label{bii2}
\end{eqnarray}
We say that a multiset $\la{\mathbb K},\sigma_{{\mathbb K}}\ra$ is empty and 
denote it by $\la{\mathbb K},\widehat{0}\;\ra$ if the following equality holds,
\begin{eqnarray}
\sigma_{{\mathbb K}}(\omega)=0\;\;\mbox{for}\;\;\forall\;\;\omega\in{\mathbb K}
\;,\;\;\;\mbox{i.e.}\;\;\;\sigma_{{\mathbb K}}=\widehat{0}:\;{\mathbb K}
\mapsto\{0\}\;,\label{z7}
\end{eqnarray}
and put for empty set ${\mathbb K}=\emptyset$ by definition $\sigma_{\emptyset}
:\;\emptyset\mapsto\{0\}$, i.e. $\sigma_{\emptyset}=\widehat{0}$ and a multiset 
$\la\emptyset,\sigma_{\emptyset}\ra$ is empty.

\noindent
Following \cite{ges95}, denote by $\bigvee$ the {\em join} operation of two 
multisets $\la{\mathbb K}_1,\sigma_{{\mathbb K}_1}\ra$ and $\la{\mathbb K}_2,
\sigma_{{\mathbb K}_2}\ra$,
\begin{eqnarray}
\la{\mathbb K}_1,\sigma_{{\mathbb K}_1}\ra\;\bigvee\;\la{\mathbb K}_2,\sigma_{
{\mathbb K}_2}\ra=\la{\mathbb K}_1\cup {\mathbb K}_2,\sigma_{{\mathbb K}_1\cup
{\mathbb K}_2}\ra\;,\;\;\;\;\;\sigma_{{\mathbb K}_1\cup{\mathbb K}_2}(\omega)=
\sigma_{{\mathbb K}_1}(\omega)+\sigma_{{\mathbb K}_2}(\omega)\;,\label{bi3}
\end{eqnarray}
so that by (\ref{bii1}) and (\ref{bi3}) we get $\la{\mathbb K}_i,\sigma_{
{\mathbb K}_i}\ra\subseteq\la{\mathbb K}_1,\sigma_{{\mathbb K}_1}\ra\;\bigvee\;
\la{\mathbb K}_2,\sigma_{{\mathbb K}_2}\ra$, $i=1,2$ and
\begin{eqnarray}
\#\left[\la{\mathbb K}_1,\sigma_{{\mathbb K}_1}\ra\;\bigvee\;\la{\mathbb K}_2,
\sigma_{{\mathbb K}_2}\ra\right]=\#\la{\mathbb K}_1,\sigma_{{\mathbb K}_1}\ra+
\#\la{\mathbb K}_2,\sigma_{{\mathbb K}_2}\ra\;.\nonumber
\end{eqnarray}
If a multiset $\la{\mathbb K}_2,\widehat{0}\;\ra$ is empty, then by (\ref{z7}) 
and (\ref{bi3}) we have $\la{\mathbb K}_1,\sigma_{{\mathbb K}_1}\ra=\la{\mathbb 
K}_1,\sigma_{{\mathbb K}_1}\ra\bigvee\la{\mathbb K}_2,\widehat{0}\;\ra=\la
{\mathbb K}_2,\widehat{0}\;\ra\bigvee\la{\mathbb K}_1,\sigma_{{\mathbb K}_1}
\ra$, that encompasses also the case ${\mathbb K}_2=\emptyset$. By the $\bigvee$
operation a multiset $\la{\mathbb K},\sigma_{{\mathbb K}}\ra$ and an element 
$\left(\omega,\sigma_{{\mathbb K}_1\cup {\mathbb K}_2}(\omega)\right)$ can be 
represented as follows,
\begin{eqnarray}
\la{\mathbb K},\sigma_{{\mathbb K}}\ra=\bigvee_{\omega\in{\mathbb K}}\left(
\omega,\sigma_{{\mathbb K}}(\omega)\right)\;,\;\;\;\;\left(\omega,\sigma_{
{\mathbb K}_1}(\omega)+\sigma_{{\mathbb K}_2}(\omega)\right)=\left(\omega,
\sigma_{{\mathbb K}_1}(\omega)\right)\bigvee\left(\omega,\sigma_{{\mathbb K}_2}
(\omega)\right)\;.\label{z10}
\end{eqnarray}
By consequence of (\ref{bi3}) the $\bigvee$ - operation satisfies the 
commutative and associative laws,
\begin{eqnarray}
\la{\mathbb K}_1,\sigma_{{\mathbb K}_1}\ra\;\bigvee\;\la{\mathbb K}_2,\sigma_{
{\mathbb K}_2}\ra=\la{\mathbb K}_2,\sigma_{{\mathbb K}_2}\ra\;\bigvee\;\la{
\mathbb K}_1,\sigma_{{\mathbb K}_1}\ra,\;\;\;\;\;\;\;\;\;\;\;\;\;\;\;\;\;\;\;\;
\;\;\;\;\;\;\;\;\label{jo1a}\\
\left[\la{\mathbb K}_1,\sigma_{{\mathbb K}_1}\ra\bigvee\la{\mathbb K}_2,\sigma
_{{\mathbb K}_2}\ra\right]\bigvee\la{\mathbb K}_3,\sigma_{{\mathbb K}_3}\ra=
\la{\mathbb K}_1,\sigma_{{\mathbb K}_1}\ra\bigvee\left[\la{\mathbb K}_2,\sigma
_{{\mathbb K}_2}\ra\bigvee\la{\mathbb K}_3,\sigma_{{\mathbb K}_3}\ra\right]=
\bigvee_{j=1}^3\la{\mathbb K}_j,\sigma_{{\mathbb K}_j}\ra.\nonumber
\end{eqnarray}
Let two multisets $\la{\mathbb K}_1,\sigma_{{\mathbb K}_1}\ra$ and $\la{\mathbb 
K}_2,\sigma_{{\mathbb K}_2}\ra$ be given such that ${\mathbb K}_1,{\mathbb K}
_2\subseteq{\mathbb N}$. i.e. if $\omega\in {\mathbb K}_1$ and $\xi\in{\mathbb
K}_2$ then $\omega+\xi\in {\mathbb N}$. Denote by ${\widehat\bigoplus}$ their 
{\em join sumset} operation and defene it as follows,
\begin{eqnarray}
&&\la{\mathbb K}_1,\sigma_{{\mathbb K}_1}\ra{\widehat\bigoplus}\la{\mathbb K}_2,
\sigma_{{\mathbb K}_2}\ra=:\bigvee_{\omega\in{\mathbb K}_1\atop\xi\in {\mathbb 
K}_2}\left[(\omega,\sigma_{{\mathbb K}_1}(\omega))\oplus(\xi,\sigma_{{\mathbb K}
_2}(\xi))\right]=\bigvee_{\omega\in{\mathbb K}_1\atop\xi\in {\mathbb K}_2}\left(
\omega+\xi,\;\sigma_{{\mathbb K}_1}(\omega)\cdot\sigma_{{\mathbb K}_2}(\xi)
\right)\;,\;\;\;\;\;\;\;\label{z1}\\
&&\#\left[\la{\mathbb K}_1,\sigma_{{\mathbb K}_1}\ra{\widehat\bigoplus}\la
{\mathbb K}_2,\sigma_{{\mathbb K}_2}\ra\right]=\#\la{\mathbb K}_1,\sigma_{
{\mathbb K}_1}\ra\cdot\#\la{\mathbb K}_2,\sigma_{{\mathbb K}_2}\ra\;,\label{z1a}
\end{eqnarray}
where $\oplus$ is a usual sumset operation which was already used in 
(\ref{bet3a}) and (\ref{tr22}). By (\ref{z1}) the ${\widehat\bigoplus}$ 
operation satisfies the commutative law,
\begin{eqnarray}   
\la{\mathbb K}_1,\sigma_{{\mathbb K}_1}\ra{\widehat\bigoplus}\la{\mathbb K}_2,
\sigma_{{\mathbb K}_2}\ra=\la{\mathbb K}_2,\sigma_{{\mathbb K}_2}\ra{\widehat
\bigoplus}\la{\mathbb K}_1,\sigma_{{\mathbb K}_1}\ra\;.\label{z5}
\end{eqnarray}
Distributive law of the ${\widehat\bigoplus}$ - operation holds over the 
$\bigvee$ - operation,
\begin{eqnarray}
\left[\la{\mathbb K}_1,\sigma_{{\mathbb K}_1}\ra\bigvee\la{\mathbb K}_2,\sigma
_{{\mathbb K}_2}\ra\right]{\widehat\bigoplus}\la{\mathbb K}_3,\sigma_{{\mathbb 
K}_3}\ra=\left[\la{\mathbb K}_1,\sigma_{{\mathbb K}_1}\ra{\widehat\bigoplus}\la
{\mathbb K}_3,\sigma_{{\mathbb K}_3}\ra\right]\bigvee\left[\la{\mathbb K}_2,
\sigma_{{\mathbb K}_2}\ra{\widehat\bigoplus}\la{\mathbb K}_3,\sigma_{{\mathbb 
K}_3}\ra\right].\label{z2}
\end{eqnarray}
We prove (\ref{z2}) making use of (\ref{bi3}), (\ref{z10}) and (\ref{z1}) and 
start with its l.h.s. in the form,
\begin{eqnarray}
\bigvee_{\omega\in{\mathbb K}_1\cup{\mathbb K}_2,\;\xi\in {\mathbb K}_3}
\left(\omega,\sigma_{{\mathbb K}_1\cup{\mathbb K}_2}(\omega)\right)\oplus
\left(\xi,\sigma_{{\mathbb K}_3}(\xi)\right)=\bigvee_{\omega\in{\mathbb K}_1
\cup{\mathbb K}_2,\;\xi\in {\mathbb K}_3}\left(\omega+\xi,\;\left[\sigma_{
{\mathbb K}_1}(\omega)+\sigma_{{\mathbb K}_2}(\omega)\right]\cdot \sigma_{
{\mathbb K}_3}(\xi)\right)\nonumber\\
=\left[\bigvee_{\omega\in{\mathbb K}_1,\;\xi\in {\mathbb K}_3}\left(\omega+\xi,
\;\sigma_{{\mathbb K}_1}(\omega)\cdot \sigma_{{\mathbb K}_3}(\xi)\right)\right]
\;\bigvee\;\left[\bigvee_{\omega\in{\mathbb K}_2,\;\xi\in {\mathbb K}_3}
\left(\omega+\xi,\;\sigma_{{\mathbb K}_2}(\omega)\cdot \sigma_{{\mathbb K}_3}
(\xi)\right)\right]\nonumber\\
=\left[\la{\mathbb K}_1,\sigma_{{\mathbb K}_1}\ra{\widehat\bigoplus}\la{\mathbb 
K}_3,\sigma_{{\mathbb K}_3}\ra\right]\bigvee\left[\la{\mathbb K}_2,\sigma_{
{\mathbb K}_2}\ra{\widehat\bigoplus}\la{\mathbb K}_3,\sigma_{{\mathbb K}_3}\ra
\right].\nonumber
\end{eqnarray}
Define an intersection of two multisets and denote it by $\bigwedge$,
\begin{eqnarray}
\la{\mathbb K}_1,\sigma_{{\mathbb K}_1}\ra\;\bigwedge\;\la{\mathbb K}_2,\;\sigma
_{{\mathbb K}_2}\ra=\la{\mathbb K}_1\cap{\mathbb K}_2,\sigma_{{\mathbb K}_1
\bigcap{\mathbb K}_2}\ra,\;\;\;\;\;\sigma_{{\mathbb K}_1\bigcap{\mathbb K}_2}
(\omega)=\min\left\{\sigma_{{\mathbb K}_1}(\omega),\sigma_{{\mathbb K}_2}
(\omega)\right\}\;,\label{bi3a}
\end{eqnarray}
so that by (\ref{bii1}) and (\ref{bi3a}) a following containment holds, 
\begin{eqnarray}
\la{\mathbb K}_1,\sigma_{{\mathbb K}_1}\ra\bigwedge\la{\mathbb K}_2,\sigma_{
{\mathbb K}_2}\ra\sqsubseteq\la{\mathbb K}_i,\sigma_{{\mathbb K}_i}\ra\;,\;\;
\;\;i=1,2\;.\label{z19}
\end{eqnarray} 
Let two multisets $\la{\mathbb K}_1,\sigma_{{\mathbb K}_1}\ra$ and $\la{\mathbb 
K}_2,\sigma_{{\mathbb K}_2}\ra$ be given such that $\la{\mathbb K}_1,\sigma_{
{\mathbb K}_1}\ra\sqsubseteq\la{\mathbb K}_2,\sigma_{{\mathbb K}_2}\ra$ in 
accordance with (\ref{bii1}). Define their set difference  $\setminus$ as 
follows,
\begin{eqnarray}
\la{\mathbb K}_2,\sigma_{{\mathbb K}_2}\ra\setminus\la{\mathbb K}_1,\sigma_{
{\mathbb K}_1}\ra=\la{\mathbb K}_2,\sigma_{{\mathbb K}_2\setminus{\mathbb K}_1}
\ra\;,\;\;\;\;\;\sigma_{{\mathbb K}_2\setminus{\mathbb K}_1}(\omega)=\sigma_{
{\mathbb K}_2}(\omega)-\sigma_{{\mathbb K}_1}(\omega)\;,\;\;\;\mbox{for}\;\;
\forall\;\;\omega\in {\mathbb K}_2\;.\label{bii3}
\end{eqnarray}
By (\ref{bi3}) and (\ref{bii3}) it follows, 
\begin{eqnarray}
\left[\la{\mathbb K}_1,\sigma_{{\mathbb K}_1}\ra\bigvee\la{\mathbb K}_2,\sigma_
{{\mathbb K}_2}\ra\right]\setminus\la{\mathbb K}_2,\sigma_{{\mathbb K}_2}\ra=
\la{\mathbb K}_1,\sigma_{{\mathbb K}_1}\ra\;,\;\;\;\la{\mathbb K},\sigma_{
{\mathbb K}}\ra\setminus\la{\mathbb K},\sigma_{{\mathbb K}}\ra=\la{\mathbb K},
\widehat{0}\;\ra\;.\label{z30}
\end{eqnarray}
We prove three Lemmas before going to the main Theorem \ref{the4} on multiset 
equalities.
\begin{lemma}\label{lem3}
Let three multisets be given such that $\la{\mathbb K}_1,\sigma_{{\mathbb K}_1}
\ra\sqsubseteq\la{\mathbb K}_2,\sigma_{{\mathbb K}_2}\ra\bigvee\la{\mathbb K}_3,
\sigma_{{\mathbb K}_3}\ra$. Then
\begin{eqnarray}
\la{\mathbb K}_1,\sigma_{{\mathbb K}_1}\ra\sqsubseteq\left[\la{\mathbb K}_1,
\sigma_{{\mathbb K}_1}\ra\bigwedge\la{\mathbb K}_2,\sigma_{{\mathbb K}_2}\ra
\right]\bigvee\left[\la{\mathbb K}_1,\sigma_{{\mathbb K}_1}\ra\bigwedge\la
{\mathbb K}_3,\sigma_{{\mathbb K}_3}\ra\right]\;.\label{z11}
\end{eqnarray}
\end{lemma}
{\sf Proof} $\;\;\;$Start with given containment and make worth of (\ref{bii2})
and (\ref{bi3})
\begin{eqnarray}
\la{\mathbb K}_1,\sigma_{{\mathbb K}_1}\ra\sqsubseteq\la{\mathbb K}_2,\sigma_{
{\mathbb K}_2}\ra\bigvee\la{\mathbb K}_3,\sigma_{{\mathbb K}_3}\ra=\la{\mathbb 
K}_2\cup{\mathbb K}_3,\sigma_{{\mathbb K}_2\cup{\mathbb K}_3}\ra\;\;\rightarrow
\;\;\left\{\begin{array}{l}{\mathbb K}_1\subseteq{\mathbb K}_2\cup{\mathbb K}_3
\;,\\\sigma_{{\mathbb K}_1}(\omega)\leq\sigma_{{\mathbb K}_2}(\omega)+\sigma_{
{\mathbb K}_3}(\omega)\end{array}\right..\label{z12}
\end{eqnarray}
However, the last containment, ${\mathbb K}_1\subseteq{\mathbb K}_2\cup{\mathbb 
K}_3$ implies a set identity ${\mathbb K}_1=\left[{\mathbb K}_1\cap{\mathbb K}_
2\right]\cup\left[{\mathbb K}_1\cap{\mathbb K}_3\right]$.

We write inequality in (\ref{z12}) in more details in three different regions 
of the set ${\mathbb K}_1$,
\begin{eqnarray}
\mbox{for}\;\;\forall\;\;\omega\in {\mathbb K}_1\;:\;\;\;\;\;
\sigma_{{\mathbb K}_1}(\omega)\leq\left\{\begin{array}{lll}
\sigma_{{\mathbb K}_2}(\omega)\;,&\;\mbox{if}\;&\omega\not\in 
{\mathbb K}_1\cap{\mathbb K}_3\;,\\
\sigma_{{\mathbb K}_2}(\omega)+\sigma_{{\mathbb K}_3}(\omega)\;,&\;\mbox{if}\;&
\omega\in{\mathbb K}_1\cap{\mathbb K}_2\cap{\mathbb K}_3\;,\\
\sigma_{{\mathbb K}_3}(\omega)\;,&\;\mbox{if}\;&\omega\not\in
{\mathbb K}_1\cap{\mathbb K}_2\;.\end{array}\right.\label{z14}
\end{eqnarray}
Apply (\ref{bi3}), (\ref{bi3a}) and identity ${\mathbb K}_1=\left[{\mathbb K}_1
\cap{\mathbb K}_2\right]\cup\left[{\mathbb K}_1\cap{\mathbb K}_3\right]$ to the 
r.h.s. of (\ref{z11}),
\begin{eqnarray}
\la{\mathbb K}_1\cap{\mathbb K}_2,\sigma_{{\mathbb K}_1\cap{\mathbb K}_2}\ra
\bigvee\la{\mathbb K}_1\cap{\mathbb K}_3,\sigma_{{\mathbb K}_1\cap{\mathbb K}_3
}\ra=\la\;\left[{\mathbb K}_1\cap{\mathbb K}_2\right]\cup\left[{\mathbb K}_1\cap
{\mathbb K}_3\right],\sigma_{{\mathbb K}_1\cap{\mathbb K}_2}+\sigma_{{\mathbb 
K}_1\cap{\mathbb K}_3}\;\ra\nonumber\\
=\la{\mathbb K}_1,\mu_{{\mathbb K}_1,{\mathbb K}_2}+\mu_{{\mathbb K}_1,{\mathbb 
K}_3}\ra\;,\;\;\;\mbox{where}\;\;\;\mu_{{\mathbb K}_i,{\mathbb K}_j}=
\min\left\{\sigma_{{\mathbb K}_i}(\omega),\sigma_{{\mathbb K}_j}(\omega)\right\}
.\label{z16}
\end{eqnarray}
Consider $\mu_{{\mathbb K}_1,{\mathbb K}_2}+\mu_{{\mathbb K}_1,{\mathbb K}
_3}$ in three different regions of the set ${\mathbb K}_1$: 1) $\omega\not\in
{\mathbb K}_1\cap{\mathbb K}_3$, 2) $\omega\not\in{\mathbb K}_1\cap{\mathbb K}_
2$ and 3) $\omega\in{\mathbb K}_1\cap{\mathbb K}_2\cap{\mathbb K}_3$. By 
(\ref{z8}) and (\ref{z14}) we have in the first two regions, 
\begin{eqnarray}
&&\mu_{{\mathbb K}_1,{\mathbb K}_2}=\left\{\begin{array}{lll}
\sigma_{{\mathbb K}_1}(\omega)&\mbox{if}&\omega\not\in
{\mathbb K}_1\cap{\mathbb K}_3\\0&\mbox{if}&\omega\not\in
{\mathbb K}_1\cap{\mathbb K}_2\end{array}\right.\;,\;\;\;\;\;
\mu_{{\mathbb K}_1,{\mathbb K}_3}=\left\{\begin{array}{lll}
0&\mbox{if}&\omega\not\in{\mathbb K}_1\cap{\mathbb K}_3\\
\sigma_{{\mathbb K}_1}(\omega)&\mbox{if}\;&\omega\not\in
{\mathbb K}_1\cap{\mathbb K}_2\end{array}\right.\;,\nonumber\\
&&\mu_{{\mathbb K}_1,{\mathbb K}_2}+\mu_{{\mathbb K}_1,{\mathbb K}_3}=\sigma_{
{\mathbb K}_1}(\omega)\;\;\;\mbox{if}\;\;\;\omega\not\in{\mathbb K}_1\cap
{\mathbb K}_2\;\;\;\mbox{or}\;\;\;\omega\not\in{\mathbb K}_1\cap{\mathbb K}_3\;
.\label{z15}
\end{eqnarray}
The 3rd region, $\omega\in{\mathbb K}_1\cap{\mathbb K}_2\cap{\mathbb K}_3$, 
requires more accurate operation. By inequality (\ref{z12}) we have three 
options: $\sigma_{{\mathbb K}_1}(\omega)\leq\sigma_{{\mathbb K}_2}(\omega)$ and 
$\sigma_{{\mathbb K}_1}(\omega)\leq\sigma_{{\mathbb K}_3}(\omega)$, or  
$\sigma_{{\mathbb K}_1}(\omega)\leq\sigma_{{\mathbb K}_2}(\omega)$ and $
\sigma_{{\mathbb K}_1}(\omega)\geq\sigma_{{\mathbb K}_3}(\omega)$ or $\sigma_{
{\mathbb K}_1}(\omega)\geq\sigma_{{\mathbb K}_2}(\omega)$ and $\sigma_{{\mathbb 
K}_1}(\omega)\leq\sigma_{{\mathbb K}_3}(\omega)$, so
\begin{eqnarray}
\mbox{for}\;\;\forall\;\;\omega\in {\mathbb K}_1\cap{\mathbb K}_2\cap{\mathbb K}
_3\;:\;\;\;
\mu_{{\mathbb K}_1,{\mathbb K}_2}+\mu_{{\mathbb K}_1,{\mathbb K}_3}=\left\{
\begin{array}{lll}2\sigma_{{\mathbb K}_1}(\omega)&\mbox{if}&
\sigma_{{\mathbb K}_1}\leq\sigma_{{\mathbb K}_2},\sigma_{{\mathbb K}_1}\leq 
\sigma_{{\mathbb K}_3},\\
\sigma_{{\mathbb K}_1}(\omega)+\sigma_{{\mathbb K}_3}(\omega)&\mbox{if}&\sigma_
{{\mathbb K}_1}\leq\sigma_{{\mathbb K}_2},\sigma_{{\mathbb K}_1}\geq\sigma_{
{\mathbb K}_3},\\
\sigma_{{\mathbb K}_1}(\omega)+\sigma_{{\mathbb K}_2}(\omega)&\mbox{if}&\sigma_
{{\mathbb K}_1}\geq\sigma_{{\mathbb K}_2},\sigma_{{\mathbb K}_1}\leq\sigma_{
{\mathbb K}_3}.\end{array}\right.\nonumber
\end{eqnarray}
Thus, keeping in mind the mid line in (\ref{z14}) we can summarize the last 
equalities and (\ref{z15}),
\begin{eqnarray}
\mbox{for}\;\;\forall\;\;\omega\in {\mathbb K}_1\;:\;\;\;\;\;
\mu_{{\mathbb K}_1,{\mathbb K}_2}+\mu_{{\mathbb K}_1,{\mathbb K}_3}\left\{
\begin{array}{lll}=\sigma_{{\mathbb K}_1}(\omega)\;,&\;\mbox{if}\;&\omega\not
\in{\mathbb K}_1\cap{\mathbb K}_3\;,\\
\geq\sigma_{{\mathbb K}_1}(\omega)\;,&\;\mbox{if}\;&\omega\in{\mathbb K}_1\cap
{\mathbb K}_2\cap{\mathbb K}_3\;,\\
=\sigma_{{\mathbb K}_1}(\omega)\;,&\;\mbox{if}\;&\omega\not\in{\mathbb K}_1
\cap{\mathbb K}_2\;.\end{array}\right.\label{z17}
\end{eqnarray}
Substituting (\ref{z17}) into (\ref{z16}) and comparing the obtained multiset 
with $\la{\mathbb K}_1,\sigma_{{\mathbb K}_1}\ra$ we arrive at (\ref{z11}) that 
finishes proof of Lemma.$\;\;\;\;\;\;\Box$

Note that a containment ${\mathbb K}_1\subseteq{\mathbb K}_2\cup{\mathbb K}_3$ 
implies a set equality ${\mathbb K}_1=\left[{\mathbb K}_1\cap{\mathbb K}_2
\right]\cup\left[{\mathbb K}_1\cap{\mathbb K}_3\right]$, but according to Lemma 
\ref{lem3} such implication cannot be extended onto multisets.

Before going to the next Lemma show that 
\begin{eqnarray}
\mbox{If}\;\;\;\la{\mathbb K}_1,\sigma_{{\mathbb K}_1}\ra\sqsubseteq\la{\mathbb 
K}_2,\sigma_{{\mathbb K}_2}\ra\sqsubseteq\la{\mathbb K}_1,\sigma_{{\mathbb K}_1
}\ra\;\;\;\mbox{then}\;\;\;\la{\mathbb K}_1,\sigma_{{\mathbb K}_1}\ra=\la
{\mathbb K}_2,\sigma_{{\mathbb K}_2}\ra\;.\label{z31}
\end{eqnarray}
Indeed, according to (\ref{bii1}) the double containment of multisets implies 
the double containment of sets, ${\mathbb K}_1\subseteq{\mathbb K}_2\subseteq
{\mathbb K}_1$, and two nonstrict inequalities, $\sigma_{{\mathbb K}_1}(\omega)
\leq\sigma_{{\mathbb K}_2}(\omega)\leq\sigma_{{\mathbb K}_1}(\omega)$, that 
gives (\ref{z31}).
\begin{lemma}\label{lem4}
Let three multisets be given such that $\la{\mathbb K}_1,\sigma_{{\mathbb K}_1}
\ra\sqsubseteq\la{\mathbb K}_2,\sigma_{{\mathbb K}_2}\ra\bigvee\la{\mathbb K}_3,
\sigma_{{\mathbb K}_3}\ra$. If a multiset $\la{\mathbb K}_1,\sigma_{{\mathbb K}
_1}\ra\bigwedge\la{\mathbb K}_2,\sigma_{{\mathbb K}_2}\ra$ is empty then
\begin{eqnarray}
\la{\mathbb K}_1,\sigma_{{\mathbb K}_1}\ra=\la{\mathbb K}_1,\sigma_{{\mathbb K}
_1}\ra\bigwedge\la{\mathbb K}_3,\sigma_{{\mathbb K}_3}\ra\;.\label{z18}
\end{eqnarray}
\end{lemma}
{\sf Proof} $\;\;\;$First, since a multiset $\la{\mathbb K}_1,\sigma_{{\mathbb 
K}_1}\ra\bigwedge\la{\mathbb K}_2,\sigma_{{\mathbb K}_2}\ra$ is empty, then 
Lemma \ref{lem3} and summation law with empty multiset imply a containment 
$\la{\mathbb K}_1,\sigma_{{\mathbb K}_1}\ra\sqsubseteq\la{\mathbb K}_1,\sigma_{
{\mathbb K}_1}\ra\bigwedge\la{\mathbb K}_3,\sigma_{{\mathbb K}_3}\ra$. However, 
by (\ref{z19}) we have an opposite containment, $\la{\mathbb K}_1,\sigma_{
{\mathbb K}_1}\ra\bigwedge\la{\mathbb K}_3,\sigma_{{\mathbb K}_3}\ra\sqsubseteq
\la{\mathbb K}_1,\sigma_{{\mathbb K}_1}\ra$. Combining together both 
containments and applying (\ref{z31}) we arrive at (\ref{z18}).$\;\;\;\;\;\;
\Box$
\begin{lemma}\label{lem5}
Let two multisets $\la{\mathbb K}_i,\sigma_{{\mathbb K}_i}\ra$, $i=1,2$ be given
such that their intersection $\la{\mathbb K}_{12},\sigma_{{\mathbb K}_{12}}\ra=
\la{\mathbb K}_1,\sigma_{{\mathbb K}_1}\ra\bigwedge{\mathbb K}_2,\sigma_{
{\mathbb K}_2}\ra$ is not empty. Then
\begin{eqnarray}
\la{\mathbb K}_i,\sigma_{{\mathbb K}_i}\ra=\left[\la{\mathbb K}_i,\sigma_{
{\mathbb K}_i}\ra\setminus\la{\mathbb K}_{12},\sigma_{{\mathbb K}_{12}}\ra
\right]\bigvee\la{\mathbb K}_{12},\sigma_{{\mathbb K}_{12}}\ra\;.\label{z21}
\end{eqnarray}
and a following multiset $\left[\la{\mathbb K}_1,\sigma_{{\mathbb K}_1}\ra
\setminus\la{\mathbb K}_{12},\sigma_{{\mathbb K}_{12}}\ra\right]\bigwedge
\left[\la{\mathbb K}_2,\sigma_{{\mathbb K}_2}\ra\setminus\la{\mathbb K}_{12},
\sigma_{{\mathbb K}_{12}}\ra\right]$ is empty,
\end{lemma}
{\sf Proof} $\;\;\;$Calculate the r.h.s. of (\ref{z21}) in accordance with 
definitions of the join '$\bigvee$' (\ref{bi3}), intersection '$\bigwedge$'
(\ref{bi3a}) and set difference '$\setminus$' (\ref{bii3}) operations for 
multisets. For $i=1$ we have
\begin{eqnarray}
\la{\mathbb K}_1,\sigma_{{\mathbb K}_1}\ra\setminus\la{\mathbb K}_{12},\sigma_
{{\mathbb K}_{12}}\ra=\la{\mathbb K}_1,\sigma_{{\mathbb K}_1}\ra\setminus\la
{\mathbb K}_1\cap{\mathbb K}_2,\sigma_{{\mathbb K}_1\cap{\mathbb K}_2}\ra=
\la{\mathbb K}_1,\sigma_{{\mathbb K}_1}-\sigma_{{\mathbb K}_1\cap{\mathbb K}_2}
\ra\;.\label{z23}
\end{eqnarray}
Thus, by consequence of (\ref{z23}) we get
\begin{eqnarray}
\left[\la{\mathbb K}_1,\sigma_{{\mathbb K}_1}\ra\setminus\la{\mathbb K}_{12},
\sigma_{{\mathbb K}_{12}}\ra\right]\bigvee\la{\mathbb K}_{12},\sigma_{{\mathbb 
K}_{12}}\ra=\la{\mathbb K}_1\cup{\mathbb K}_{12},\sigma_{{\mathbb K}_1}-\sigma_
{{\mathbb K}_1\cap {\mathbb K}_2}+\sigma_{{\mathbb K}_1\cap{\mathbb K}_2}\ra=
\la{\mathbb K}_1,\sigma_{{\mathbb K}_1}\ra\;.\nonumber
\end{eqnarray}
The proof for $i=2$ is similar. Thus, the 1st part of Lemma is proven. 
As for the 2nd part of Lemma, consider a multiset $\la{\mathbb K}_{12},\sigma_{
{\mathbb K}_{12}}\ra$ according to (\ref{bi3a}),
\begin{eqnarray}
\la{\mathbb K}_{12},\sigma_{{\mathbb K}_{12}}\ra=\la{\mathbb K}_1\cap{\mathbb K}
_2,\sigma_{{\mathbb K}_1\bigcap{\mathbb K}_2}\ra\;,\;\;\;\;\sigma_{{\mathbb K}_1
\bigcap{\mathbb K}_2}(\omega)=\left\{\begin{array}{lll}\sigma_{{\mathbb K}_1}
(\omega),&\mbox{if}\;&\sigma_{{\mathbb K}_1}(\omega)\leq\sigma_{{\mathbb K}_2}
(\omega),\\\sigma_{{\mathbb K}_2}(\omega),&\mbox{if}&\sigma_{{\mathbb K}_1}
(\omega)\geq\sigma_{{\mathbb K}_2}(\omega).\end{array}\right.\label{z24}
\end{eqnarray}
Then, by (\ref{z23}) and (\ref{z24}) we have
\begin{eqnarray}
\la{\mathbb K}_1,\sigma_{{\mathbb K}_1}\ra\setminus\la{\mathbb K}_{12},\sigma_
{{\mathbb K}_{12}}\ra=\la{\mathbb K}_1,\sigma_{{\mathbb K}_1}-\sigma_{{\mathbb 
K}_1\cap{\mathbb K}_2}\ra\;,\;\;\;\la{\mathbb K}_2,\sigma_{{\mathbb K}_2}\ra
\setminus\la{\mathbb K}_{12},\sigma_{{\mathbb K}_{12}}\ra=\la{\mathbb K}_2,
\sigma_{{\mathbb K}_2}-\sigma_{{\mathbb K}_1\cap{\mathbb K}_2}\ra,\label{z25}
\end{eqnarray}
\begin{eqnarray}
&&\mbox{If}\;\;\;\sigma_{{\mathbb K}_1}(\omega)\leq\sigma_{{\mathbb K}_2}
(\omega)\;\;\;\mbox{then}\;\;\;\sigma_{{\mathbb K}_1}-\sigma_{{\mathbb K}_1
\cap{\mathbb K}_2}=0\;,\;\;\sigma_{{\mathbb K}_2}-\sigma_{{\mathbb K}_1\cap
{\mathbb K}_2}=\sigma_{{\mathbb K}_2}(\omega)-\sigma_{{\mathbb K}_1}(\omega)\;,
\nonumber\\
&&\mbox{If}\;\;\;\sigma_{{\mathbb K}_1}(\omega)\geq\sigma_{{\mathbb K}_2}
(\omega)\;\;\;\mbox{then}\;\;\;\sigma_{{\mathbb K}_1}-\sigma_{{\mathbb K}_1\cap
{\mathbb K}_2}=\sigma_{{\mathbb K}_1}(\omega)-\sigma_{{\mathbb K}_2}(\omega)\;,
\;\;\sigma_{{\mathbb K}_2}-\sigma_{{\mathbb K}_1\cap{\mathbb K}_2}=0\;.
\label{z26}
\end{eqnarray}
Combining (\ref{z25}) and (\ref{z26}) with the $\bigwedge$ operation 
(\ref{bi3a}) we can calculate an intersection
\begin{eqnarray}
\left[\la{\mathbb K}_1,\sigma_{{\mathbb K}_1}\ra\setminus\la{\mathbb K}_{12},  
\sigma_{{\mathbb K}_{12}}\ra\right]\bigwedge\left[\la{\mathbb K}_2,\sigma_{
{\mathbb K}_2}\ra\setminus\la{\mathbb K}_{12},\sigma_{{\mathbb K}_{12}}\ra
\right]=\la{\mathbb K}_1\cap{\mathbb K}_2,\widehat{0}\;\ra\;,\nonumber
\end{eqnarray}
that means in accordance with (\ref{z7}) that multiset $\left[\la{\mathbb K}_1,
\sigma_{{\mathbb K}_1}\ra\setminus\la{\mathbb K}_{12},\sigma_{{\mathbb K}_{12}}
\ra\right]\bigwedge\left[\la{\mathbb K}_2,\sigma_{{\mathbb K}_2}\ra\setminus
\la{\mathbb K}_{12},\sigma_{{\mathbb K}_{12}}\ra\right]$ is empty. Thus, our 
Lemma is proven completely.$\;\;\;\;\;\;\Box$

In this paper we study two multisets $\la{\mathbb M},\sigma_{{\mathbb M}}\ra$ 
and $\la{\mathbb X},\sigma_{{\mathbb X}}\ra$ described at the end of section 
\ref{s3}. As for $\la{\mathbb M},\sigma_{{\mathbb M}}\ra$, the multiplicity 
$\sigma_{{\mathbb M}}(\omega)$ accounts for the number of elements $\omega=d_{
i_1}+\ldots+d_{i_k}+F_r$, $\omega\in{\mathbb M}$, $1\leq k\leq m-1$, $1\leq r
\leq \gamma\left({\bf d}^m\right)$, of equal values. As for $\la{\mathbb X},
\sigma_{{\mathbb X}}\ra$, the multiplicity $\sigma_{{\mathbb X}}(\xi)$ accounts 
for the number of elements $\xi=C_{j,i},{\overline C}_{j,i}$, $\xi\in{\mathbb 
X}$, $1\leq i\leq m-1$, $1\leq j\leq\beta_i\left({\bf d}^m\right)$, of equal 
values that comes by much more difficult way making use of the Hilbert syzygy 
theorem \cite{eis05}.

In fact, in this paper we deal mainly with cardinalities of entire multisets $\la{\mathbb M},
\sigma_{{\mathbb M}}\ra$ and $\la{\mathbb X},\sigma_{{\mathbb X}}\ra$ or of 
their submultisets. Therefore, for the sake of brevity we will often reduce the 
designations for multisets $\la{\mathbb M},\sigma_{{\mathbb M}}\ra={\mathfrak 
M}$ and $\la{\mathbb X},\sigma_{{\mathbb X}}\ra={\mathfrak X}$, skipping the 
underlying sets ${\mathbb M}$, ${\mathbb X}$ and the mapping functions $\sigma_
{{\mathbb M}}$, $\sigma_{{\mathbb X}}$. For example, we shall write (\ref{z2}) 
as follows, $\left[{\mathfrak K}_1\bigvee{\mathfrak K}_2\right]{\widehat
\bigoplus}\;{\mathfrak K}_3=\left[{\mathfrak K}_1\;{\widehat\bigoplus}\;
{\mathfrak K}_3\right]\bigvee\left[{\mathfrak K}_2\;{\widehat\bigoplus}\;
{\mathfrak K}_3\right]$, where ${\mathfrak K}_i=\la{\mathbb K}_i,\sigma_{
{\mathbb K}_i}\ra$. We hope that such reduction will not mislead the readers.
\subsection{Multisets and equation (\ref{tr34b})}\label{s41}
After recasting the terms of Eq. (\ref{tr34b}) in such a way that in every its 
l.h.s. and r.h.s. would remain only positive terms, all degrees in power terms 
can be arranged in 4 multisets ${\mathfrak M}_1$, ${\mathfrak M}_2$ and 
${\mathfrak X}_1$, ${\mathfrak X}_2$,
\begin{eqnarray}
\sum_{\omega_1\in{\mathbb M}_1\left({\bf d}^m\right)}z^{\omega_1}+
\sum_{\xi_1\in{\mathbb X}_1\left({\bf d}^m\right)}z^{\xi_1}=
\sum_{\omega_2\in{\mathbb M}_2\left({\bf d}^m\right)}z^{\omega_2}+
\sum_{\xi_2\in{\mathbb X}_2\left({\bf d}^m\right)}z^{\xi_2}\;.\label{z40}
\end{eqnarray}
Two of these multisets, ${\mathfrak M}_1=\la{\mathbb M}_1,\sigma_{{\mathbb M}_1
}\ra$ and ${\mathfrak X}_1=\la{\mathbb X}_1,\sigma_{{\mathbb X}_1}\ra$, are 
distributed in the l.h.s. of Eq. (\ref{tr34b}) while the other two, ${\mathfrak 
M}_2=\la{\mathbb M}_2,\sigma_{{\mathbb M}_2}\ra$ and ${\mathfrak X}_2=
\la{\mathbb X}_2,\sigma_{{\mathbb X}_2}\ra$, are distributed in the r.h.s., of 
Eq. (\ref{tr34b}). The sets ${\mathbb M}_1$ and ${\mathbb M}_2$ are the sets of 
partial sums $\omega_{1,2}=d_{i_1}+\ldots+d_{i_k}+F_r$ of gaps $F_j$ and 
generators $d_i$. Both sets ${\mathbb X}_1$ and ${\mathbb X}_2$ are the sets of 
degrees $\xi_{1,2}=C_{j,i},\;{\overline C}_{j,i}$ defined in (\ref{ar10a}) and 
(\ref{z35}). In view of definition (\ref{bi3}) of the operation $\bigvee$ the 
multiset equality associated with Eq. (\ref{z40}) reads,
\begin{eqnarray}
{\mathfrak M}_1\bigvee{\mathfrak X}_1={\mathfrak M}_2\bigvee{\mathfrak X}_2\;.
\label{jo1c}
\end{eqnarray}
Denote by ${\mathfrak \not O}$ the empty multiset and prove the following 
theorem on multiset equalities.
\begin{theorem}\label{the4}Let two finite multisets ${\mathfrak M}_1$, 
${\mathfrak M}_2$ of integers and two finite multisets ${\mathfrak X}_1$, 
${\mathfrak X}_2$ of indeterminate elements be given such that
\begin{eqnarray}
{\mathfrak M}_1\bigwedge {\mathfrak M}_2={\mathfrak M}_{12}\neq{\mathfrak \not 
O}\;,\;\;\;\;{\mathfrak X}_1\bigwedge{\mathfrak X}_2={\mathfrak X}_{12}\neq
{\mathfrak \not O}\;,\label{jo1b}
\end{eqnarray}
and let a multiset equality (\ref{jo1c}) be given. Then the following hold
\begin{eqnarray}
{\mathfrak X}_1\setminus{\mathfrak X}_{12}={\mathfrak M}_2\setminus{\mathfrak 
M}_{12}\;,\;\;\;\;\;\;\;{\mathfrak X}_2\setminus {\mathfrak X}_{12}={\mathfrak 
M}_1\setminus{\mathfrak M}_{12}\;.\label{jo1d}
\end{eqnarray}
\end{theorem}
{\sf Proof} $\;\;\;$By consequence of (\ref{z21}) we can represent the multisets
${\mathfrak M}_1$, ${\mathfrak M}_2$ and ${\mathfrak X}_1$, ${\mathfrak X}_2$ 
as follows,
\begin{eqnarray}
{\mathfrak M}_i=\left[{\mathfrak M}_i\setminus{\mathfrak M}_{12}\right]\bigvee
{\mathfrak M}_{12}\;,\;\;\;\;\;\;\;{\mathfrak X}_i=\left[{\mathfrak X}_i
\setminus{\mathfrak X}_{12}\right]\bigvee{\mathfrak X}_{12}\;,\;\;\;\;i=1,2\;,
\label{jo1e}
\end{eqnarray}
and substitute (\ref{jo1e}) into (\ref{jo1c}),
\begin{eqnarray}
\left[\left[{\mathfrak M}_1\setminus {\mathfrak M}_{12}\right]\bigvee{\mathfrak 
M}_{12}\right]\bigvee\left[\left[{\mathfrak X}_1\setminus{\mathfrak X}_{12}
\right]\bigvee{\mathfrak X}_{12}\right]=\left[\left[{\mathfrak M}_2\setminus
{\mathfrak M}_{12}\right]\bigvee{\mathfrak M}_{12}\right]\bigvee\left[\left[
{\mathfrak X}_2\setminus{\mathfrak X}_{12}\right]\bigvee{\mathfrak X}_{12}
\right].\nonumber
\end{eqnarray}
Making use of commutative and associative laws (\ref{jo1a}) rewrite the last 
equation as follows,
\begin{eqnarray}
\left[\left[{\mathfrak M}_1\setminus{\mathfrak M}_{12}\right]\bigvee\left[
{\mathfrak X}_1\setminus {\mathfrak X}_{12}\right]\right]\bigvee\left[{\mathfrak
X}_{12}\bigvee{\mathfrak M}_{12}\right]=\left[\left[{\mathfrak M}_2\setminus
{\mathfrak M}_{12}\right\}\bigvee\left[{\mathfrak X}_2\setminus{\mathfrak X}_{
12}\right]\right]\bigvee\left[{\mathfrak X}_{12}\bigvee{\mathfrak M}_{12}
\right].\label{jo1u}
\end{eqnarray}
According to (\ref{z30}) take a complement of multiset ${\mathfrak X}_{12}
\bigvee{\mathfrak M}_{12}$ in the l.h.s. and r.h.s. of (\ref{jo1u}),
\begin{eqnarray}
\left[{\mathfrak M}_1\setminus {\mathfrak M}_{12}\right]\bigvee
\left[{\mathfrak X}_1\setminus {\mathfrak X}_{12}\right]=
\left[{\mathfrak M}_2\setminus {\mathfrak M}_{12}\right]\bigvee
\left[{\mathfrak X}_2\setminus {\mathfrak X}_{12}\right]\;.\label{jo1f}
\end{eqnarray}
However, in accordance with the 2nd part of Lemma \ref{lem5} two following 
pairs of multisets are disjoined,
\begin{eqnarray}
\left[{\mathfrak M}_1\setminus{\mathfrak M}_{12}\right]\bigwedge\left[{\mathfrak
M}_2\setminus{\mathfrak M}_{12}\right]={\mathfrak \not O}\;,\;\;\;\;\left[
{\mathfrak X}_1\setminus {\mathfrak X}_{12}\right]\bigwedge\left[{\mathfrak X}_
2\setminus {\mathfrak X}_{12}\right]={\mathfrak \not O}\;.\label{jo1g}
\end{eqnarray}
Comparing multiset equality (\ref{jo1f}) supplied with conditions (\ref{jo1g}) 
and Lemma \ref{lem4} we arrive at (\ref{jo1d}).$\;\;\;\;\;\;\Box$
\subsection{Multisets ${\mathfrak M}_1$, ${\mathfrak M}_2$ and their 
intersection ${\mathfrak M}_{12}$}\label{s42}
In this section we give a detailed description of the multisets ${\mathfrak M}_
1$ and ${\mathfrak M}_2$ which were introduced in section \ref{s41}. The 
multisets ${\mathfrak X}_1$ and ${\mathfrak X}_2$ will be constructed in 
sections \ref{s5} ($m=2n$) and \ref{s6} ($m=2n+1$). 

Consider the 1st line in Eq. (\ref{tr34b}) which is the only giving rise to 
multisets ${\mathfrak M}_1$ and ${\mathfrak M}_2$,
\begin{eqnarray}
\Lambda\left({\bf d}^m;z\right)=\sum_{r=1}^{\gamma\left({\bf d}^m\right)}
\left[\sum_{j=1}^mz^{d_j+F_r}-\sum_{j>k=1}^mz^{d_j+d_k+F_r}+\ldots+(-1)^m
\sum_{j=1}^mz^{\Sigma_m-d_j+F_r}\right].\;\label{z37}
\end{eqnarray}
Making use of definitions (\ref{bi3}) and (\ref{z1}) of the operations $\bigvee$
and ${\widehat\bigoplus}$, construct a sequence of multisets 
\begin{eqnarray}
{\mathfrak L}_i\left({\bf d}^m\right)={\mathfrak D}_i\left({\bf d}^m\right)\;
{\widehat\bigoplus}\;\Delta_{{\cal H}}\left({\bf d}^m\right),\;i=\leq m-1\;,\;\;
\;\mbox{where}\;\;\;\;{\mathfrak D}_i\left({\bf d}^m\right)=\la{\mathbb D}_i
\left({\bf d}^m\right),\sigma_{{\mathbb D}_i}\ra\;.\label{jo3}
\end{eqnarray}
The underlying sets ${\mathbb D}_i\left({\bf d}^m\right)$ of degrees $x_j$ and 
the mapping functions $\sigma_{{\mathbb D}_i}\left(x_j\right)$ are given by
\begin{eqnarray}
{\mathbb D}_1\left({\bf d}^m\right)=\bigcup_{i=1}^m\{d_i\},\;\;\;{\mathbb D}_2
\left({\bf d}^m\right)=\bigcup_{i_1>i_2=1}^m\{d_{i_1}+d_{i_2}\},\;\;\;{\mathbb 
D}_3\left({\bf d}^m\right)=\bigcup_{i_1>i_2>i_3=1}^m\{d_{i_1}+d_{i_2}+d_{i_3}\},
\;\ldots\nonumber\\
\sigma_{{\mathbb D}_1}\left(d_i\right)=1\;,\;\;\;\sigma_{{\mathbb D}_2}\left(
d_{i_1}+d_{i_2}\right)=1\;,\;\;\;\sigma_{{\mathbb D}_3}\left(d_{i_1}+d_{i_2}
+d_{i_3}\right)=1\;,\;\;\ldots\;.\;\;\;\;\;\;\label{z63}
\end{eqnarray}
Due to (\ref{z63}) we have $\#{\mathfrak D}_i\left({\bf d}^m\right)={m\choose 
i}$ and according to (\ref{z1a}) the entire cardinality $\#{\mathfrak L}_i
\left({\bf d}^m\right)$ reads,
\begin{eqnarray}
\#{\mathfrak L}_i\left({\bf d}^m\right)=\gamma\left({\bf d}^m\right)\cdot 
\#{\mathfrak D}_i\left({\bf d}^m\right)\;,\;\;\;\mbox{and by definition}\;\;\;
\#{\mathfrak L}_0\left({\bf d}^m\right)=0\;.\label{z36}
\end{eqnarray}
Continue to compose multisets and construct two other auxiliary multisets,
\begin{eqnarray}
{\mathfrak L}_o\left({\bf d}^m\right)={\mathfrak L}_1\left({\bf d}^m\right)
\;\bigvee\;{\mathfrak L}_3\left({\bf d}^m\right)\;\bigvee\;\ldots\;\bigvee
\;{\mathfrak L}_{2\left\lfloor\frac{m}{2}\right\rfloor-1}\left({\bf d}^m\right)
\;,\;\;\;\;\;\;{\mathfrak L}_o\left({\bf d}^m\right)=\la{\mathbb L}_o
\left({\bf d}^m\right),\sigma_{{\mathbb L}_o}\ra\;,\;\label{tr34y}\\
{\mathfrak L}_e\left({\bf d}^m\right)={\mathfrak L}_2\left({\bf d}^m\right)
\;\bigvee\;{\mathfrak L}_4\left({\bf d}^m\right)\;\bigvee\;\ldots\;\bigvee
\;{\mathfrak L}_{2\left\lfloor\frac{m-1}{2}\right\rfloor}\left({\bf d}^m\right)
\;,\;\;\;\;\;\;{\mathfrak L}_e\left({\bf d}^m\right)=\la{\mathbb L}_e
\left({\bf d}^m\right),\sigma_{{\mathbb L}_e}\ra\;.\;\;\nonumber
\end{eqnarray}
In (\ref{tr34y}) subscripts '$o$' and '$e$' stand for the odd and even numbers 
of summands $d_j$ in the elements of ${\mathfrak L}_k\left({\bf d}^m\right)$, 
respectively. By consequence of distributive property (\ref{z2}) of ${\widehat
\bigoplus}$ over $\bigvee$ we have,
\begin{eqnarray}
&&{\mathfrak L}_o\left({\bf d}^m\right)={\mathfrak D}_o\left({\bf d}^m\right)\;
{\widehat\bigoplus}\;\Delta_{{\cal H}}\left({\bf d}^m\right)\;,\;\;\;\;\;
{\mathfrak D}_o\left({\bf d}^m\right)={\mathfrak D}_1\left({\bf d}^m\right)\;
\bigvee\;\ldots\;\bigvee\;{\mathfrak D}_{2\left\lfloor\frac{m}{2}\right
\rfloor-1}\left({\bf d}^m\right)\;,\;\;\;\label{lp21}\\
&&{\mathfrak L}_e\left({\bf d}^m\right)={\mathfrak D}_e\left({\bf d}^m\right)\;
{\widehat\bigoplus}\;\Delta_{{\cal H}}\left({\bf d}^m\right)\;,\;\;\;\;\;
{\mathfrak D}_e\left({\bf d}^m\right)={\mathfrak D}_2\left({\bf d}^m\right)\;
\bigvee\;\ldots\;\bigvee\;{\mathfrak D}_{2\left\lfloor\frac{m-1}{2}\right
\rfloor}\left({\bf d}^m\right)\;.\;\;\;\nonumber
\end{eqnarray}
Represent the 1st line in Eq. (\ref{tr34b}) (see (\ref{z37})) in terms of the 
multisets elements $\lambda_o$ and $\lambda_e$,
\begin{eqnarray}
\Lambda\left({\bf d}^m;z\right)=\sum_{\lambda_o\in{\mathbb L}_o\left({\bf d}^m
\right)}z^{\lambda_o}-\sum_{\lambda_e\in{\mathbb L}_e\left({\bf d}^m\right)}
z^{\lambda_e}\;,\label{z38}
\end{eqnarray}
where two underlying sets ${\mathbb L}_o\left({\bf d}^m\right)$ and ${\mathbb L}
_e\left({\bf d}^m\right)$ were introduced in (\ref{tr34y}). 

Keeping in mind the definition of the multisets ${\mathfrak M}_1$ and 
${\mathfrak M}_2$ given in section \ref{s41} (after Eq. (\ref{z40})) and 
comparing the difference $\Lambda\left({\bf d}^m;z\right)$ in (\ref{z38}) coming
from (\ref{z37}) with recasted Eq. (\ref{z40}) we conclude that ${\mathfrak L}_
o\left({\bf d}^m\right)$, ${\mathfrak L}_e\left({\bf d}^m\right)$ and 
${\mathfrak L}_e\left({\bf d}^m\right)\bigwedge{\mathfrak L}_o\left({\bf d}^m
\right)$ are exactly the multisets ${\mathfrak M}_1$, ${\mathfrak M}_2$ and 
${\mathfrak M}_{12}$ appeared in (\ref{jo1c}) and Theorem \ref{the4} and 
equipped with dependence on ${\bf d}^m$,
\begin{eqnarray}
{\mathfrak M}_1\left({\bf d}^m\right)\equiv {\mathfrak L}_o\left({\bf d}^m
\right)\;,\;\;\;\;\;{\mathfrak M}_2\left({\bf d}^m\right)\equiv {\mathfrak L}_
e\left({\bf d}^m\right)\;,\;\;\;\;\;{\mathfrak M}_{12}\left({\bf d}^m\right)
\equiv{\mathfrak L}_e\left({\bf d}^m\right)\bigwedge{\mathfrak L}_o\left({\bf 
d}^m\right)\;.\label{lp22}
\end{eqnarray}
\subsection{Multiset ${\mathfrak D}_e\left({\bf d}^m\right)\bigwedge{\mathfrak 
D}_o\left({\bf d}^m\right)$ for small edim}\label{s43}
In the case of pseudosymmetric semigroup (see section \ref{s21}) the expressions
for ${\mathfrak M}_{12}\left({\bf d}^m\right)$ and its cardinality are 
simplified considerably. Indeed, since $\#\Delta_{{\cal H}}\left({\bf d}^m
\right)=1$ then by (\ref{lp21}) we get
\begin{eqnarray}
{\mathfrak L}_o\left({\bf d}^m\right)={\mathfrak D}_o\left({\bf d}^m\right)
\oplus\left\{F\left({\bf d}^m\right)/2\right\}\;,\;\;\;\;
{\mathfrak L}_e\left({\bf d}^m\right)={\mathfrak D}_e\left({\bf d}^m\right)
\oplus\left\{F\left({\bf d}^m\right)/2\right\}\;.\nonumber
\end{eqnarray}
Define a new multiset, ${\mathfrak D}_{eo}\left({\bf d}^m\right)={\mathfrak D}_
e\left({\bf d}^m\right)\bigwedge{\mathfrak D}_o\left({\bf d}^m\right)$ and 
represent ${\mathfrak M}_{12}\left({\bf d}^m\right)$ according to (\ref{lp22}) 
\begin{eqnarray}
{\mathfrak M}_{12}\left({\bf d}^m\right)={\mathfrak D}_{eo}\left({\bf d}^m
\right)\oplus\left\{1/2\;F({\bf d}^m)\right\}\;,\;\;\;\;\mbox{and}\;\;\;\;
\#{\mathfrak M}_{12}\left({\bf d}^m\right)=\#{\mathfrak D}_{eo}\left({\bf d}^m
\right)\;.\nonumber
\end{eqnarray}
In the sequel we address the following questions: how big can be the cardinality
$\#{\mathfrak D}_{eo}\left({\bf d}^m\right)$ and when it does vanish. For 
generic tuple ${\bf d}^m$ these questions are addressed to the additive number 
theory. Here we give answer for small edim, $m=3,4,5$, and return to arbitrary 
edim elsewhere.
\begin{remark}\label{rem1}Making use of representation (\ref{z10}) for multisets
${\mathfrak D}_k\left({\bf d}^m\right)$ and unique occurrence (\ref{z63}) of any
sum $\omega_k=d_{i_1}+d_{i_2}+\ldots +d_{i_k}$, $1\leq k<m$, of generators $d_j$
therein we skip (in this section) a unity $\sigma_{{\mathbb D}_k}(\omega_k)=1$ 
in representation (\ref{z10}). E.g. we shall write, 

$\;\;\;\;\;\;\;\;\;\;\;\;{\mathfrak D}_2\left({\bf d}^4\right)=\left\{d_1+d_2,
d_1+d_3,d_1+d_4,d_2+d_3,d_2+d_4,d_3+d_4\right\}$.
\end{remark}
\begin{proposition}\label{pro2}Two multisets ${\mathfrak D}_{eo}\left({\bf d}^3
\right)$ and ${\mathfrak D}_{eo}\left({\bf d}^4\right)$ are empty.
\end{proposition}
{\sf Proof} $\;\;\;$First, according to (\ref{lp21}) write ${\mathfrak D}_o
\left({\bf d}^3\right)={\mathfrak D}_1\left({\bf d}^3\right)$ and ${\mathfrak D}
_e\left({\bf d}^3\right)={\mathfrak D}_2\left({\bf d}^3\right)$, so
\begin{eqnarray}
{\mathfrak D}_o\left({\bf d}^3\right)=\left\{d_1,d_2,d_3\right\}\;,\;\;\;\;
{\mathfrak D}_e\left({\bf d}^3\right)=\left\{d_1+d_2,d_2+d_3,d_3+d_1\right\}\;.
\nonumber
\end{eqnarray}
Thus, ${\mathfrak D}_{eo}\left({\bf d}^3\right)={\mathfrak \not O}$, otherwise 
the minimality of the generating set ${\bf d}^3$ would be broken. Next, 
according to (\ref{lp21}) write ${\mathfrak D}_o\left({\bf d}^4\right)=
{\mathfrak D}_1\left({\bf d}^4\right)\bigvee{\mathfrak D}_3\left({\bf d}^4
\right)$ and ${\mathfrak D}_e\left({\bf d}^4\right)={\mathfrak D}_2\left({\bf 
d}^4\right)$, so
\begin{eqnarray}
&&{\mathfrak D}_o\left({\bf d}^4\right)=\left\{d_1,d_2,d_3,d_4,d_1+d_2+d_3,d_2+
d_3+d_4,d_3+d_4+d_1,d_1+d_3+d_4\right\}\;,\nonumber\\
&&{\mathfrak D}_e\left({\bf d}^4\right)=\left\{d_1+d_2,d_1+d_3,d_1+d_4,d_2+d_3,
d_2+d_4,d_3+d_4\right\}\;.\nonumber
\end{eqnarray}
Thus, ${\mathfrak D}_{eo}\left({\bf d}^4\right)={\mathfrak \not O}$ by the same 
reason of minimality of the generating set ${\bf d}^4$.$\;\;\;\;\;\;\Box$
\begin{proposition}\label{pro3}Let a numerical semigroup ${\sf S}\left({\bf d}^5
\right)$ be given. Then $\#{\mathfrak D}_{eo}\left({\bf d}^5\right)\leq 1$.
\end{proposition}
{\sf Proof} $\;\;\;$Let a tuple ${\bf d}^5$ be given, then according to
(\ref{lp21}) write
\begin{eqnarray}
{\mathfrak D}_o\left({\bf d}^5\right)={\mathfrak D}_1\left({\bf d}^5\right)
\bigvee{\mathfrak D}_3\left({\bf d}^5\right)\;,\;\;\;\;
{\mathfrak D}_e\left({\bf d}^5\right)={\mathfrak D}_2\left({\bf d}^5\right)
\bigvee{\mathfrak D}_4\left({\bf d}^5\right)\;.\nonumber
\end{eqnarray}
It is easy to verify that
\begin{eqnarray}
{\mathfrak D}_1\left({\bf d}^5\right)\bigwedge{\mathfrak D}_2\left({\bf d}^5
\right)={\mathfrak D}_1\left({\bf d}^5\right)\bigwedge{\mathfrak D}_4\left(
{\bf d}^5\right)={\mathfrak D}_3\left({\bf d}^5\right)\bigwedge{\mathfrak D}_4
\left({\bf d}^5\right)={\mathfrak \not O}\;,\nonumber
\end{eqnarray}
by the reason of minimality of the generating set ${\bf d}^5$. There is left a 
multiset ${\mathfrak D}_2\left({\bf d}^5\right)\bigwedge{\mathfrak D}_3\left(
{\bf d}^5\right)$ which we'll study. Write all ${5\choose 2}{5-2\choose 3}=10$ 
admissible relations and check their compatibility,
\begin{eqnarray}
\left.\begin{array}{rrr}1)\;\;d_1+d_2=d_3+d_4+d_5\;,&2)\;\;d_1+d_3=d_2+d_4+d_5
\;,&3)\;\;d_1+d_4=d_2+d_3+d_5\;,\\4)\;\;d_1+d_5=d_2+d_3+d_4\;,&5)\;\;d_2+d_3= 
d_1+d_4+d_5\;,&6)\;\;d_2+d_4=d_1+d_3+d_5\;,\\7)\;\;d_2+d_5=d_1+d_3+d_4\;,&8)\;
\;d_3+d_4=d_1+d_2+d_5\;,&9)\;\;d_3+d_5=d_1+d_2+d_4\;,\\10)\;\;d_4+d_5=d_1+d_2+
d_3\;.& & \end{array}\right.\label{w6w}
\end{eqnarray}
Since the ordering in the set $\{1,2,3,4,5\}$ is arbitrary it is sufficient to 
check the compatability of the 1st equality in (\ref{w6w}) with the other nine. 
By inspection of compatibility of all 9 pairs we conclude that it contradicts 
the minimality of the generating set $\{d_1,d_2,d_3,d_4,d_5\}$, e. g.
\begin{eqnarray}
\left\{\begin{array}{r}1)\;\;d_1+d_2=d_3+d_4+d_5\\2)\;\;d_1+d_3=d_2+d_4+d_5   
\end{array}\right.\rightarrow\;d_2=d_3\;,\;\;\;\;\left\{\begin{array}{r}1)\;\; 
d_1+d_2=d_3+d_4+d_5\\9)\;\;d_3+d_5=d_1+d_2+d_4\end{array}\right.\rightarrow   
\;d_5=0\;,\;\;\;\;\mbox{etc}\;.\nonumber
\end{eqnarray}
Thus, there exists at most one admissible relation and Proposition is proven.
$\;\;\;\;\;\;\Box$
\section{Almost Symmetric Semigroups ${\sf S}\left({\bf d}^{2n}\right)$ , 
$n\geq 2$}\label{s5}
In (\ref{jo1c}) we have defined a multiset equality associated with Eq. 
(\ref{tr34b}) and based on two multisets ${\mathfrak M}_1$, ${\mathfrak M}_2$ 
of given gaps $F_j$ and generators $d_i$, and two multisets ${\mathfrak X}_1$, 
${\mathfrak X}_2$ of degrees $C_{j,i}$ and ${\overline C}_{j,i}$. The first two 
multisets ${\mathfrak M}_1$ and ${\mathfrak M}_2$ were constructed explicitly in
(\ref{lp22}). In this section we construct the other two multisets ${\mathfrak 
X}_1$ and ${\mathfrak X}_2$ providing their consistence with Eqs. (\ref{tr34b}) 
and (\ref{jo1c}).

An interchange of signs of the terms in Eq. (\ref{tr34b}) and factor $(-1)^{m-
1}$ make our analysis not easy, this can be seen for edim of distinct parities, 
$m=2n$ and $m=2n+1$, where the multisets ${\mathfrak X}_1$ and ${\mathfrak X}_2$
are composed in different ways. Therefore we consider two cases of even and odd 
edim separately, and start with for $m=2n$. Substituting into the 1st line of 
Eq. (\ref{tr34b}) its representation given in (\ref{z38}) and (\ref{lp22}), 
write the whole Eq. (\ref{tr34b}) in the form which is similar to (\ref{z40}),
\begin{eqnarray}
\sum_{\omega_1\in{\mathbb M}_1\left({\bf d}^{2n}\right)}z^{\omega_1}+\sum_{q=
1}^{n-1}\sum_{j=1}^{\beta_{2q-1}\left({\bf d}^{2n}\right)}\left[z^{C_{j,2q-1}}+
z^{{\overline C}_{j,2q-1}}\right]=\sum_{\omega_2\in{\mathbb M}_2\left({\bf d}
^{2n}\right)}z^{\omega_2}+\sum_{q=1}^{n-1}\sum_{j=1}^{\beta_{2q}\left({\bf d}^{
2n}\right)}\left[z^{C_{j,2q}}+z^{{\overline C}_{j,2q}}\right]\label{tr34d}
\end{eqnarray}
By comparison of Eqs. (\ref{tr34d}) and (\ref{z40}) we'll find the multisets 
${\mathfrak X}_1\left({\bf d}^{2n}\right)$ and ${\mathfrak X}_2\left({\bf d}^
{2n}\right)$. 

Consider the last sums in the l.h.s. and r.h.s. of Eq. (\ref{tr34d}) and 
construct two auxiliary multisets ${\mathfrak B}_i\left({\bf d}^{2n}\right)$ 
and ${\overline{\mathfrak B}}_i\left({\bf d}^{2n}\right)$ of syzygies degrees
which have a standard representations (see section \ref{s4}) through the sets 
${\mathbb B}_i\left({\bf d}^{2n}\right)$ and ${\overline{\mathbb B}}_i\left(
{\bf d}^{2n}\right)$ defined in (\ref{ar10a}) and (\ref{z35}), 
\begin{eqnarray}
{\mathfrak B}_i\left({\bf d}^{2n}\right)=\la{\mathbb B}_i\left({\bf d}^{2n}
\right),\sigma_{{\mathbb B}_i\left({\bf d}^{2n}\right)}\ra\;,\;\;\;\;
{\overline{\mathfrak B}}_i\left({\bf d}^{2n}\right)=\la{\overline{\mathbb B}}_i
\left({\bf d}^{2n}\right),\sigma_{{\overline{\mathbb B}}_i\left({\bf d}^{2n}
\right)}\ra\;.\label{tr34p}
\end{eqnarray}
By definitions (\ref{ar10a}) and (\ref{z35}) of the sets ${\mathbb B}_i\left(
{\bf d}^{2n}\right)$ and ${\overline{\mathbb B}}_i\left({\bf d}^{2n}\right)$
we have $\sigma_{{\mathbb B}_i\left({\bf d}^{2n}\right)}\left(C_{j,i}\right)=
\sigma_{{\overline{\mathbb B}}_i\left({\bf d}^{2n}\right)}\left({\overline C}_
{j,i}\right)$ that together with (\ref{tr34p}) leads to $\#{\mathfrak B}_i
\left({\bf d}^{2n}\right)=\#{\overline{\mathfrak B}}_i\left({\bf d}^{2n}\right)
=\beta_i\left({\bf d}^{2n}\right)$. Define the following multisets,
\begin{eqnarray}
{\mathfrak B}_o\left({\bf d}^{2n}\right)=\bigvee_{i=1}^{n-1}{\mathfrak B}_{2i-1}
\left({\bf d}^{2n}\right)\;,\;\;\;\;{\mathfrak B}_e\left({\bf d}^{2n}\right)=
\bigvee_{i=1}^{n-1}{\mathfrak B}_{2i}\left({\bf d}^{2n}\right)\;,\label{tr34x}\\
{\overline{\mathfrak B}}_o\left({\bf d}^{2n}\right)=\bigvee_{i=1}^{n-1}
{\overline{\mathfrak B}}_{2i-1}\left({\bf d}^{2n}\right)\;,\;\;\;\;
{\overline{\mathfrak B}}_e\left({\bf d}^{2n}\right)=\bigvee_{i=1}^{n-1}
{\overline{\mathfrak B}}_{2i}\left({\bf d}^{2n}\right)\;,\nonumber
\end{eqnarray}
where subscripts '$o$' and '$e$' stand for the odd $q=2i-1$ and even $q=2i$ 
indices, respectively, of summands ${\mathfrak B}_q\left({\bf d}^{2n}\right)$ 
and ${\overline{\mathfrak B}}_q\left({\bf d}^{2n}\right)$.

By comparison of Eqs. (\ref{tr34d}) and (\ref{z40}) we can define the multisets 
${\mathfrak X}_1\left({\bf d}^{2n}\right)$, ${\mathfrak X}_2\left({\bf d}^{2n}
\right)$ and their intersection ${\mathfrak X}_{12}\left({\bf d}^{2n}\right)$ 
through four multisets (\ref{tr34x}) and two multiset operations $\bigvee$ and 
$\bigwedge$,
\begin{eqnarray}
{\mathfrak X}_1\left({\bf d}^{2n}\right)\equiv {\mathfrak B}_o\left({\bf d}^{
2n}\right)\;\bigvee\;{\overline{\mathfrak B}}_o\left({\bf d}^{2n}\right)\;,\;\;
\;\;\;\;\;{\mathfrak X}_2\left({\bf d}^{2n}\right)\equiv{\mathfrak B}_e\left(
{\bf d}^{2n}\right)\;\bigvee\;{\overline{\mathfrak B}}_e\left({\bf d}^{2n}
\right)\;,\label{w13}\\
{\mathfrak X}_{12}\left({\bf d}^{2n}\right)\equiv\left[{\mathfrak B}_o\left({\bf
d}^{2n}\right)\;\bigvee\;{\overline{\mathfrak B}}_o\left({\bf d}^{2n}\right)
\right]\bigwedge\left[{\mathfrak B}_e\left({\bf d}^{2n}\right)\;\bigvee\;   
{\overline{\mathfrak B}}_e\left({\bf d}^{2n}\right)\right]\;.\;\;\;\;\;
\label{lp2}
\end{eqnarray}
Substituting (\ref{lp22}) and (\ref{w13}) into multiset equality (\ref{jo1c}) 
we get
\begin{eqnarray}
{\mathfrak L}_{o}\left({\bf d}^{2n}\right)\;\bigvee\;{\mathfrak B}_{o}\left(
{\bf d}^{2n}\right)\;\bigvee\;{\overline{\mathfrak B}}_o\left({\bf d}^{2n}
\right)={\mathfrak L}_{e}\left({\bf d}^{2n}\right)\;\bigvee\;{\mathfrak B}_e
\left({\bf d}^{2n}\right)\;\bigvee\;{\overline{\mathfrak B}}_e\left({\bf d}^{
2n}\right)\;.\nonumber
\end{eqnarray}
\begin{lemma}\label{lem6}Let an almost symmetric semigroup ${\sf S}\left({\bf
d}^{2n}\right)$ be given. Then 
\begin{eqnarray}
&&\left[{\mathfrak B}_o\left({\bf d}^{2n}\right)\;\bigvee\;{\overline{\mathfrak 
B}}_o\left({\bf d}^{2n}\right)\right]\setminus{\mathfrak X}_{12}\left({\bf d}^{
2n}\right)={\mathfrak L}_e\left({\bf d}^{2n}\right)\setminus{\mathfrak M}_{12}
\left({\bf d}^{2n}\right)\;,\label{jo3a}\\
&&\left[{\mathfrak B}_e\left({\bf d}^{2n}\right)\;\bigvee\;{\overline{\mathfrak 
B}}_e\left({\bf d}^{2n}\right)\right]\setminus{\mathfrak X}_{12}\left({\bf d}^{
2n}\right)={\mathfrak L}_o\left({\bf d}^{2n}\right)\setminus{\mathfrak M}_{12}
\left({\bf d}^{2n}\right)\;.\nonumber
\end{eqnarray}
\end{lemma}
{\sf Proof} $\;\;\;$Substituting the expressions (\ref{lp22}) for multisets 
${\mathfrak M}_i\left({\bf d}^{2n}\right)$, $i=1,2$, and expressions (\ref{w13})
for multisets ${\mathfrak X}_i\left({\bf d}^{2n}\right)$, $i=1,2$, into equality
(\ref{jo1c}) we apply Theorem \ref{the4}. Thus, by consequence of (\ref{jo1d}) 
we arrive at (\ref{jo3a}).$\;\;\;\;\;\;\Box$

Lemma \ref{lem6} does not give yet explicit expressions for syzygies degrees 
$C_{j,i}$ and ${\overline C}_{j,i}$ since it is hard to differentiate them one 
from another. This requires much more powerful algebraic methods, e.g. the 
Hilbert syzygy theorem \cite{eis05}. However, Lemma \ref{lem6} leads to new 
relations for the Betti numbers. Define the following cardinalities: 
$\#{\mathfrak M}_{12}\left({\bf d}^m \right)=\ell\left({\bf d}^m\right)$, 
$\#{\mathfrak X}_{12}\left({\bf d}^m \right)=\wp\left({\bf d}^m\right)$ and 
$\delta\left({\bf d}^m\right)=\wp\left({\bf d}^m\right)-\ell\left({\bf 
d}^m\right)$ and prove Theorem.
\begin{theorem}\label{the5}Let an almost symmetric semigroup ${\sf S}\left({\bf 
d}^{2n}\right)$ be given. Then
\begin{eqnarray}
&&\beta_1\left({\bf d}^{2n}\right)+\beta_3\left({\bf d}^{2n}\right)+\ldots +
\beta_{2n-3}\left({\bf d}^{2n}\right)=\gamma\left({\bf d}^{2n}\right)\cdot
\left(4^{n-1}-1\right)+\frac1{2}\;\delta\left({\bf d}^{2n}\right)\;,
\label{jo4a}\\
&&\beta_2\left({\bf d}^{2n}\right)+\beta_4\left({\bf d}^{2n}\right)+\ldots +
\beta_{2n-2}\left({\bf d}^{2n}\right)=\gamma\left({\bf d}^{2n}\right)\cdot 
4^{n-1}+\frac1{2}\;\delta\left({\bf d}^{2n}\right)\;.\nonumber
\end{eqnarray}
\end{theorem}
{\sf Proof} $\;\;\;$By Lemma \ref{lem6} and in view of definition (\ref{bi3}) 
of the operation $\bigvee$ we get
\begin{eqnarray}
&&\#{\mathfrak B}_o\left({\bf d}^{2n}\right)+\#{\overline{\mathfrak B}}_0\left(
{\bf d}^{2n}\right)-\wp\left({\bf d}^{2n}\right)=\#{\mathfrak L}_e\left({\bf d}
^{2n}\right)-\ell\left({\bf d}^{2n}\right)\;,\nonumber\\
&&\#{\mathfrak B}_e\left({\bf d}^{2n}\right)+\#{\overline{\mathfrak B}}_e\left(
{\bf d}^{2n}\right)-\wp\left({\bf d}^{2n}\right)=\#{\mathfrak L}_{o}\left({\bf 
d}^{2n}\right)-\ell\left({\bf d}^{2n}\right)\;.
\nonumber
\end{eqnarray}
Making use of (\ref{jo3}), (\ref{tr34y}) and (\ref{tr34p}), (\ref{tr34x}), and 
inserting them into the last equations we arrive at
\begin{eqnarray}
\beta_1\left({\bf d}^{2n}\right)+\beta_3\left({\bf d}^{2n}\right)+\ldots +
\beta_{2n-3}\left({\bf d}^{2n}\right)=\frac{\gamma\left({\bf d}^{2n}\right)}{2}
\left[{2n\choose 2}+{2n\choose 4}+\ldots +{2n\choose 2n-2}\right]+
\frac{\delta\left({\bf d}^{2n}\right)}{2}\;,\nonumber\\
\beta_2\left({\bf d}^{2n}\right)+\beta_4\left({\bf d}^{2n}\right)+\ldots 
+\beta_{2n-2}\left({\bf d}^{2n}\right)=\frac{\gamma\left({\bf d}^{2n}\right)}{2}
\left[{2n\choose 1}+{2n\choose 3}+\ldots +{2n\choose 2n-1}\right]+
\frac{\delta\left({\bf d}^{2n}\right)}{2}\;.\nonumber
\end{eqnarray}
A simple algebraic exercise gives,
\begin{eqnarray}
{2n\choose 2}+{2n\choose 4}+\ldots +{2n\choose 2n-2}=2^{2n-1}-2\;,\;\;
{2n\choose 1}+{2n\choose 3}+\ldots +{2n\choose 2n-1}=2^{2n-1}\;,\nonumber
\end{eqnarray}
that bring us to (\ref{jo4a}).$\;\;\;\;\;\;\Box$ 

By consequence of (\ref{jo4a}) and the fact that the Betti numbers are
nonnegative integers it follows that  $\wp\left({\bf d}^{2n}\right)=\ell\left({
\bf d}^{2n}\right)\pmod 2$.

The case of pseudosymmetric semigroup, $\gamma\left({\bf d}^m\right)=1$, of 
embedding dimension 4 is most simple. By Proposition \ref{pro2} and Theorem 
\ref{the5} we have here,
\begin{eqnarray}
\beta_1\left({\bf d}^4\right)=3+\frac1{2}\;\wp\left({\bf d}^4\right)\;,\;\;\;\;
\beta_2\left({\bf d}^4\right)=4+\frac1{2}\;\wp\left({\bf d}^4\right)\;.
\label{z33}
\end{eqnarray}
\begin{example}$\{d_1,d_2,d_3,d_4\}=\{5,6,7,9\}\;,\;\;\beta_1=5\;,\;\;
\beta_2=6\;,\;\;\beta_3=2\;,\;\;$
\label{ex1}{\footnotesize
\begin{eqnarray}
B_3(5,6,7,9)&=&\{31,35\}\;,\;\;{\sf S}^{\prime}(5,6,7,9)=\{4,8\}\;,\;\;
\Delta_{{\cal H}}(5,6,7,9)=\{4\}\;,\;\;\;\sum_4=27\;,\nonumber\\
\Delta_{{\cal G}}(5,6,7,9)&=&\{1,2,3,8\}\;,\;\;F(5,6,7,9)=8\;,\;\;
\ell(5,6,7,9)=0\;,\;\;\wp(5,6,7,9)=4\;,\nonumber\\
Q(5,6,7,9;z)&=&1-z^{12}-z^{14}-z^{15}-z^{16}-z^{18}+z^{21}+z^{22}+z^{23}+
z^{24}+z^{25}+z^{26}-z^{31}-z^{35}\;.\nonumber
\end{eqnarray}}
\end{example}
\section{Almost Symmetric Semigroups ${\sf S}\left({\bf d}^{2n+1}\right)$ , 
$n\geq 1$}\label{s6}
Substituting into the 1st line of Eq. (\ref{tr34b}) its representation given in 
(\ref{z38}) and (\ref{lp22}), write  Eq. (\ref{tr34b}) for $m=2n+1$ as follows,
\begin{eqnarray}
\sum_{\omega_1\in{\mathbb M}_1\left({\bf d}^{2n+1}\right)}z^{\omega_1}+
\sum_{q=1}^n\sum_{j=1}^{\beta_{2q-1}\left({\bf d}^{2n+1}\right)}z^{C_{j,2q-1}}+
\sum_{q=1}^{n-1}\sum_{j=1}^{\beta_{2q}\left({\bf d}^{2n+1}\right)}z^{{\overline 
C}_{j,2q}}=\;\;\;\;\;\;\;\;\;\;\;\;\;\;\;\;\label{tr34c}\\
\sum_{\omega_2\in{\mathbb M}_2\left({\bf d}^{2n+1}\right)}z^{\omega_2}+
\sum_{q=1}^{n-1}\sum_{j=1}^{\beta_{2q}\left({\bf d}^{2n+1}\right)}z^{C_{j,2q}}+
\sum_{q=1}^n\sum_{j=1}^{\beta_{2q-1}\left({\bf d}^{2n+1}\right)}z^{{\overline 
C}_{j,2q-1}}\;.\nonumber
\end{eqnarray}
By comparison of Eqs. (\ref{tr34c}) and (\ref{z40}) we'll find the multisets
${\mathfrak X}_1\left({\bf d}^{2n+1}\right)$ and ${\mathfrak X}_2\left({\bf d}^
{2n+1}\right)$.

Consider the last sums in the l.h.s. and r.h.s. of Eq. (\ref{tr34c}) and
construct two auxiliary multisets ${\mathfrak B}_i\left({\bf d}^{2n+1}\right)$
and ${\overline{\mathfrak B}}_i\left({\bf d}^{2n+1}\right)$ of syzygies degrees
which have a standard representations (see section \ref{s4}) through the sets 
${\mathbb B}_i\left({\bf d}^{2n+1}\right)$ and ${\overline{\mathbb B}}_i\left(
{\bf d}^{2n+1}\right)$ defined in (\ref{ar10a}) and (\ref{z35}),
\begin{eqnarray}
{\mathfrak B}_i\left({\bf d}^{2n+1}\right)=\la{\mathbb B}_i\left({\bf d}^{2n+1}
\right),\sigma_{{\mathbb B}_i\left({\bf d}^{2n+1}\right)}\ra\;,\;\;\;\;
{\overline{\mathfrak B}}_i\left({\bf d}^{2n+1}\right)=\la{\overline{\mathbb B}}
_i\left({\bf d}^{2n+1}\right),\sigma_{{\overline{\mathbb B}}_i\left({\bf d}^{
2n+1}\right)}\ra\;.\label{gp0}
\end{eqnarray}
By definitions (\ref{ar10a}) and (\ref{z35}) of the underlying sets we have
$\sigma_{{\mathbb B}_i\left({\bf d}^{2n+1}\right)}\left(C_{j,i}\right)=\sigma_{
{\overline{\mathbb B}}_i\left({\bf d}^{2n+1}\right)}\left({\overline C}_{j,i}
\right)$ that together with (\ref{gp0}) leads to $\#{\mathfrak B}_i\left({\bf 
d}^{2n+1}\right)=\#{\overline{\mathfrak B}}_i\left({\bf d}^{2n+1}\right)=
\beta_i\left({\bf d}^{2n+1}\right)$. Define four other multisets,
\begin{eqnarray}
{\mathfrak B}_o\left({\bf d}^{2n+1}\right)=\bigvee_{i=1}^{n}{\mathfrak B}_{2i-1}
\left({\bf d}^{2n+1}\right)\;,\;\;\;\;{\mathfrak B}_e\left({\bf d}^{2n+1}\right)
=\bigvee_{i=1}^{n-1}{\mathfrak B}_{2i}\left({\bf d}^{2n+1}\right)\;,\nonumber\\
{\overline{\mathfrak B}}_o\left({\bf d}^{2n+1}\right)=\bigvee_{i=1}^{n}
{\overline{\mathfrak B}}_{2i-1}\left({\bf d}^{2n+1}\right)\;,\;\;\;\;
{\overline{\mathfrak B}}_e\left({\bf d}^{2n+1}\right)=\bigvee_{i=1}^{n-1}
{\overline{\mathfrak B}}_{2i}\left({\bf d}^{2n+1}\right)\;,\label{gp3}
\end{eqnarray}
where subscripts '$o$' and '$e$' stand for the odd $q=2i-1$ and even $q=2i$
indices, respectively, of summands ${\mathfrak B}_q\left({\bf d}^{2n+1}\right)$ 
and ${\overline{\mathfrak B}}_q\left({\bf d}^{2n+1}\right)$.

By comparison of Eqs. (\ref{tr34c}) and (\ref{z40}) we can define the multisets 
${\mathfrak X}_1\left({\bf d}^{2n+1}\right)$, ${\mathfrak X}_2\left({\bf d}^{2n
+1}\right)$ and their intersection ${\mathfrak X}_{12}\left({\bf d}^{2n+1}
\right)$ through four multisets (\ref{gp3}) and two multiset operations 
$\bigvee$ and $\bigwedge$,
\begin{eqnarray}
{\mathfrak X}_1\left({\bf d}^{2n+1}\right)\equiv{\mathfrak B}_o\left({\bf d}^{
2n+1}\right)\;\bigvee\;{\overline{\mathfrak B}}_e\left({\bf d}^{2n+1}\right)\;,
\;\;\;\;\;\;\;{\mathfrak X}_2\left({\bf d}^{2n+1}\right)\equiv{\mathfrak B}_e
\left({\bf d}^{2n+1}\right)\;\bigvee\;{\overline{\mathfrak B}}_o\left({\bf d}^{2n+1}
\right)\;,\label{w17}\\
{\mathfrak X}_{12}\left({\bf d}^{2n+1}\right)\equiv\left[{\mathfrak B}_o\left( 
{\bf d}^{2n+1}\right)\;\bigvee\;{\overline{\mathfrak B}}_e\left({\bf d}^{2n+1}
\right)\right]\bigwedge\left[{\mathfrak B}_e\left({\bf d}^{2n+1}\right)\;\bigvee
\;{\overline{\mathfrak B}}_o\left({\bf d}^{2n+1}\right)\right].\;\;\;\;\;
\label{gp5}
\end{eqnarray}
Substituting (\ref{lp22}) and (\ref{w17}) into multiset equality (\ref{jo1c})
we get
\begin{eqnarray}
{\mathfrak L}_o\left({\bf d}^{2n+1}\right)\;\bigvee\;{\mathfrak B}_o\left({\bf d
}^{2n+1}\right)\;\bigvee\;{\overline{\mathfrak B}}_e\left({\bf d}^{2n+1}\right)=
{\mathfrak L}_e\left({\bf d}^{2n+1}\right)\;\bigvee\;{\mathfrak B}_e\left({\bf d
}^{2n+1}\right)\;\bigvee\;{\overline{\mathfrak B}}_o\left({\bf d}^{2n+1}\right)
\;.\nonumber
\end{eqnarray}
\begin{lemma}\label{lem7}Let an almost symmetric semigroup ${\sf S}\left({\bf
d}^{2n+1}\right)$ be given. Then
\begin{eqnarray}
&&\left[{\mathfrak B}_o\left({\bf d}^{2n+1}\right)\;\bigvee\;{\overline
{\mathfrak B}}_e\left({\bf d}^{2n+1}\right)\right]\setminus{\mathfrak X}_{12}
\left({\bf d}^{2n+1}\right)={\mathfrak L}_e\left({\bf d}^{2n+1}\right)
\setminus{\mathfrak M}_{12}\left({\bf d}^{2n+1}\right)\;,\label{jo5a}\\
&&\left[{\mathfrak B}_e\left({\bf d}^{2n+1}\right)\;\bigvee\;{\overline
{\mathfrak B}}_o\left({\bf d}^{2n+1}\right)\right]\setminus {\mathfrak X}_{12}
\left({\bf d}^{2n+1}\right)={\mathfrak L}_o\left({\bf d}^{2n+1}\right)\setminus
{\mathfrak M}_{12}\left({\bf d}^{2n+1}\right)\;.\nonumber
\end{eqnarray}
\end{lemma}
{\sf Proof} $\;\;\;$Substituting the expressions (\ref{lp22}) for multisets
${\mathfrak M}_i\left({\bf d}^{2n+1}\right)$, $i=1,2$, and expressions 
(\ref{w17}) for multisets ${\mathfrak X}_i\left({\bf d}^{2n+1}\right)$, $i=1,2$,
into equality (\ref{jo1c}) we apply Theorem \ref{the4}. Thus, by consequence of 
(\ref{jo1d}) we arrive at (\ref{jo5a}).$\;\;\;\;\;\;\Box$
\begin{theorem}\label{the6}Let an almost symmetric semigroup ${\sf S}\left({\bf 
d}^{2n+1}\right)$ be given. Then 
\begin{eqnarray}
&&\beta_1\left({\bf d}^{2n+1}\right)+\beta_3\left({\bf d}^{2n+1}\right)+\ldots +
\beta_{2n-1}\left({\bf d}^{2n+1}\right)={\gamma\left({\bf d}^{2n+1}\right)}
\cdot 2^{2n-1}+\frac1{2}\;\delta\left({\bf d}^{2n+1}\right)+1\;,\label{jo5c}\\
&&\beta_2\left({\bf d}^{2n+1}\right)+\beta_4\left({\bf d}^{2n+1}\right)+\ldots +
\beta_{2n-2}\left({\bf d}^{2n+1}\right)={\gamma\left({\bf d}^{2n+1}\right)}\cdot
\left(2^{2n-1}-1\right)+\frac1{2}\;\delta\left({\bf d}^{2n+1}\right)-1.\;\;\;
\;\;\nonumber
\end{eqnarray}
\end{theorem}
{\sf Proof} $\;\;\;$By Lemma \ref{lem7} and definition (\ref{bi3}) of the 
operation $\bigvee$ we have
\begin{eqnarray}
&&\#{\mathfrak B}_o\left({\bf d}^{2n+1}\right)+\#{\overline{\mathfrak B}}_e
\left({\bf d}^{2n+1}\right)-\wp\left({\bf d}^{2n+1}\right)=\#{\mathfrak L}_e
\left({\bf d}^{2n+1}\right)-\ell\left({\bf d}^{2n+1}\right)\;,
\nonumber\\
&&\#{\mathfrak B}_e\left({\bf d}^{2n+1}\right)+\#{\overline{\mathfrak B}}_o
\left({\bf d}^{2n+1}\right)-\wp\left({\bf d}^{2n+1}\right)=\#{\mathfrak L}_o
\left({\bf d}^{2n+1}\right)-\ell\left({\bf d}^{2n+1}\right)\;.
\nonumber
\end{eqnarray}
Making use of (\ref{jo3}), (\ref{tr34y}) and (\ref{gp0}), (\ref{gp3}), and 
inserting them into the last equations we get
\begin{eqnarray}
\beta_1\left({\bf d}^{2n+1}\right)+\beta_2\left({\bf d}^{2n+1}\right)+\ldots+
\beta_{2n-1}\left({\bf d}^{2n+1}\right)=\gamma\left({\bf d}^{2n+1}\right)\cdot
\left(4^n-1\right)+\delta\left({\bf d}^{2n+1}\right)\;.\nonumber
\end{eqnarray}
Making sum of the last equality with (\ref{bet2}) and simplifying the 
result we arrive at (\ref{jo5c}). $\;\;\;\;\;\;\Box$

The following Example of almost symmetric semigroups ${\sf S}\left({\bf d}^5
\right)$ is  taken from \cite{bafr97}. We have calculated the Hilbert series, 
the Betti numbers and the other characteristics.
\begin{example}$\{d_1,d_2,d_3,d_4,d_5\}=\{6,7,8,10,11\}\;,\;\;   
\beta_1=9,\;\;\beta_2=17,\;\;\beta_3=12,\;\;\beta_4=3$
\label{ex2}{\footnotesize
\begin{eqnarray}
B_4(6,7,8,10,11)&=&\{46,47,51\}\;,\;\;{\sf S}^{\prime}(6,7,8,10,11)=\{4,5,9\}
\;,\;\;\Delta_{{\cal H}}(6,7,8,10,11)=\{4,5\}\;,\;\;\sum_5=42\nonumber\\
\Delta_{{\cal G}}(6,7,8,10,11)&=&\{1,2,3,9\}\;,\;\;F(6,7,8,10,11)=9\;,\;\;
\ell(6,7,8,10,11)=2\;,\;\;\wp(6,7,8,10,11)=10\;,\nonumber\\
Q(6,7,8,10,11;z)&=&1-z^{14}-z^{16}-z^{17}-2z^{18}-z^{19}-z^{20}-z^{21}-z^{22}+
z^{24}+2z^{25}+2z^{26}+2z^{27}+3z^{28}+\nonumber\\
&&3z^{29}+2z^{30}+z^{31}+z^{32}-2z^{35}-2z^{36}-2z^{37}-z^{38}-2z^{39}-2z^{40}
-z^{41}+z^{46}+z^{47}+z^{51}\nonumber
\end{eqnarray}}
\end{example}
\begin{corollary}\label{cor2}Let an almost symmetric semigroup ${\sf S}\left(
{\bf d}^m\right)$ be given. Then
\begin{eqnarray}
\delta\left({\bf d}^m\right)\leq \left[d_1-\gamma\left({\bf d}^m\right)\right]
2^{m-1}-2m\;.\label{w80}
\end{eqnarray}
\end{corollary}
{\sf Proof} $\;\;\;$We prove Corollary for even and odd edim separately. First, 
consider an almost symmetric semigroup ${\sf S}\left({\bf d}^{2n}\right)$ and 
calculate the sum of the Betti numbers $\beta_k\left({\bf d}^{2n}\right)$. 
Keeping in mind $\beta_{2n-1}\left({\bf d}^{2n}\right)=2n-1$ and making use of 
Theorem \ref{the5} we get,
\begin{eqnarray}
\sum_{k=0}^{{2n-1}}\beta_k\left({\bf d}^{2n}\right)=\gamma\left({\bf d}^{2n}
\right)2^{2n-1}+2+\delta\left({\bf d}^{2n}\right)\;.\label{w81}
\end{eqnarray}
By comparison (\ref{w81}) with (\ref{bet22}) we obtain,
\begin{eqnarray}
\delta\left({\bf d}^{2n}\right)\leq \left[d_1-\gamma\left({\bf d}^{2n}\right)
\right]2^{2n-1}-4n\;.\label{w82}
\end{eqnarray}  
Next, consider an almost symmetric semigroup ${\sf S}\left({\bf d}^{2n+1}
\right)$ and make similar calculations with help of Theorem \ref{the6},
\begin{eqnarray}
\sum_{k=0}^{{2n}}\beta_k\left({\bf d}^{2n+1}\right)=\gamma\left({\bf d}^{2n+1}
\right)4^n+2+\delta\left({\bf d}^{2n+1}\right)\;.\label{w83}
\end{eqnarray}
By comparison (\ref{w83}) with (\ref{bet22}) we obtain,
\begin{eqnarray}
\delta\left({\bf d}^{2n+1}\right)\leq \left[d_1-\gamma\left({\bf d}^{2n+1}
\right)\right]2^{2n}-2(2n+1)\;.\label{w84}
\end{eqnarray}
Combining formulas (\ref{w82}) and (\ref{w84}) we come to (\ref{w80}).
$\;\;\;\;\;\;\Box$
\subsection{Pseudosymmetric semigroup ${\sf S}\left({\bf d}^3\right)$}
\label{s61}
This case is mostly simple and makes it possible to find all syzygy degrees and 
the Frobenius number as well. Keeping in mind $\beta_1\left({\bf d}^3\right)=3$,
$t\left({\bf d}^3\right)=2$, $F_1=\frac1{2}F\left({\bf d}^3\right)$ and denoting
$C_{j,1}=e_j$, write Eq. (\ref{tr34c})
\begin{eqnarray}
\sum_{j=1}^3z^{d_j+\frac1{2}F\left({\bf d}^3\right)}+\sum_{j=1}^3z^{e_j}=
\sum_{j>k=1}^3z^{d_j+d_k+\frac1{2}F\left({\bf d}^3\right)}+\sum_{j=1}^3z^{-   
e_j+F\left({\bf d}^3\right)+\Sigma_3}\;.\nonumber
\end{eqnarray}
However, by Proposition \ref{pro2} we have $\ell\left({\bf d}^3\right)=0$, and 
therefore by the 1st equation in (\ref{jo5c}) it follows
\begin{eqnarray}
3=3+\frac1{2}\wp\left({\bf d}^3\right)\;\;\;\;\longrightarrow\;\;\;\;\wp\left(
{\bf d}^3\right)=0\;\;\;\;\longrightarrow\;\;\;\;{\mathfrak X}_{12}\left({\bf 
d}^3\right)={\mathfrak \not O}\;.\label{w31}
\end{eqnarray}
By (\ref{w31}) and Lemma \ref{lem7} the multiset equalities (\ref{jo5a}) read: 
${\mathfrak B}_1\left({\bf d}^3\right)={\mathfrak L}_2\left({\bf d}^3\right)$
and ${\overline{\mathfrak B}}_1\left({\bf d}^3\right)={\mathfrak L}_1\left({\bf 
d}^3\right)$ that gives two following correspondences,
\begin{eqnarray}
\left\{e_i\right\}=\left\{d_j+d_k+\frac1{2}F\left({\bf d}^3\right)\right\}\;\;\;
\mbox{and}\;\;\;\left\{d_i+\frac1{2}F\left({\bf d}^3\right)\right\}=\left\{
-e_i+F\left({\bf d}^3\right)+\Sigma_3\right\}\;,\label{w33}
\end{eqnarray}
which are consistent each other. Hence, the whole numerator $Q\left({\bf d}^3;z
\right)$ in the Hilbert series reads,
\begin{eqnarray}
Q\left({\bf d}^3;z\right)=1-z^{d_1+d_2+\frac1{2}F\left({\bf d}^3\right)}-z^{d_2
+d_3+\frac1{2}F\left({\bf d}^3\right)}-z^{d_3+d_1+\frac1{2}F\left({\bf d}^3
\right)}+z^{\frac1{2}F\left({\bf d}^3\right)+\Sigma_3}+z^{F\left({\bf d}^3
\right)+\Sigma_3}.\nonumber
\end{eqnarray}
The last expression makes it possible to derive the explicit formulas for the
Frobenius number $F\left({\bf d}^3\right)$ and genus $G\left({\bf d}^3\right)$ 
of the 3D pseudosymmetric semigroups. For this purpose we'll make use of 
formulas for generic 3D nonsymmetric semigroups which were established in 
\cite{fel04}, Ch. 6,
\begin{eqnarray}
2F\left({\bf d}^3\right)=E_1+\sqrt{E_1^2-4E_2+4D_3}-2D_1\;,\;\;\;\;2G\left({\bf 
d}^3\right)=1+E_1-\frac{E_3}{D_3}-D_1\;,\label{w36}\\
\mbox{where}\;\;\;\;E_1=e_1+e_2+e_3\;,\;\;\;E_2=e_1e_2+e_2e_3+e_3e_1\;,\;\;\;
E_3=e_1e_2e_3\;,\label{w34}
\end{eqnarray}
and $D_1=\Sigma_3$, $D_2=d_1d_2+d_2d_3+d_3d_1$, $D_3=d_1d_2d_3$. Substituting 
a correspondence (\ref{w33}) into (\ref{w34}) we calculate,
\begin{eqnarray}
E_1&=&2D_1+\frac{3}{2}F\left({\bf d}^3\right)\;,\;\;\;E_2=3\left[D_1+\frac1{2}
F\left({\bf d}^3\right)\right]^2-2D_1\left[D_1+\frac1{2}F\left({\bf d}^3\right)
\right]+D_2\;,\nonumber\\
E_3&=&\left[D_1+\frac1{2}F\left({\bf d}^3\right)\right]^3-D_1\left[D_1+\frac1{
2}F\left({\bf d}^3\right)\right]^2+D_2\left[D_1+\frac1{2}F\left({\bf d}^3
\right)\right]-D_3\;.\label{w35}
\end{eqnarray}
Next, inserting (\ref{w35}) into (\ref{w36}) we get
\begin{eqnarray}
F\left({\bf d}^3\right)=-D_1+\sqrt{D_1^2+4(D_3-D_2)}\;,\;\;\;\;\;G\left({\bf d}
^3\right)=1+\frac1{2}\;F\left({\bf d}^3\right)\;.\label{w37}
\end{eqnarray}
Formulas (\ref{w37}) have been derived independently in \cite{ros05} by 
analyzing the Ap\'ery set of pseudosymmetric semigroup ${\sf S}\left({\bf d}^3
\right)$.
\section{Almost Symmetric Semigroups ${\sf S}\left({\bf d}^m\right)$ of Maximal 
edim}\label{s7}
A study of almost symmetric semigroups ${\sf S}\left({\bf d}^m\right)$ with 
maximal edim is motivated by two reasons. First, there are many known results 
\cite{heku71}, \cite{sal79}, \cite{feai07} on generic semigroups ${\sf S}\left(
{\bf d}^m_{MED}\right)$ of maximal edim (MED) that makes it reasonable to apply 
to them the statements of sections \ref{s5} and \ref{s6}. Next, Proposition 
\ref{pro1} in \cite{bafr97}, at p. 426, is followed by remark:
'{\em not any almost symmetric MED--semigroup of type $t$ and Frobenius number 
$g$ is of the type described in Propos. 11, as following example shows}',
\begin{example}$\{d_1,d_2,d_3,d_4\}=\{4,10,19,25\}\;,\;\;\beta_1=6,\;\;\beta_2=
8,\;\;\beta_3=3$\label{ex3}{\footnotesize
\begin{eqnarray}
B_4(4,10,19,25)&=&\{64,73,79\}\;,\;\;{\sf S}^{\prime}(4,10,19,25)=\{6,15,21\}
\;,\;\;\Delta_{{\cal H}}(4,10,19,25)=\{6,15\}\;,\nonumber\\
\Delta_{{\cal G}}(4,10,19,25)&=&\{1,2,3,5,7,9,11,13,17,21\}\;,\;\;
F(4,10,19,25)=21\;,\;\;G(4,10,19,25)=12\;,\nonumber\\
Q(4,10,19,25;z)&=&1-z^{20}-z^{29}-z^{35}-z^{38}+z^{39}-z^{44}+z^{45}+
z^{48}-z^{50}+2z^{54}+z^{60}+z^{63}-z^{64}+\nonumber\\
&&z^{69}-z^{73}-z^{79}\;,\;\;\;\;\;\;\;\;\;\;\ell(6,7,8,10,11)=
\wp(6,7,8,10,11)=0\;,\;\;\;\;\;\;\;\sum_4=58\;.\nonumber
\end{eqnarray}}
\end{example}
In Example \ref{ex3} we have calculated the Hilbert series, the Betti 
numbers and the other characteristics.

\noindent
Thus, there is a necessity to give the most wide description of almost symmetric
MED--semigroups.

Start with known results on the MED--semigroups \cite{heku71}, \cite{sal79} and
\cite{feai07}, Corollary 8 : 
\begin{eqnarray}
&&F\left({\bf d}^m_{MED}\right)=d_m-m\;,\;\;\;G\left({\bf d}^m_{MED}\right)=
\frac{1}{m}\sum_{k=2}^md_k-\frac{m-1}{2}\;,\;\;\;\beta_k\left({\bf d}^m_{MED}
\right)=k{m\choose k+1}\;,\;\;\;\;\;\;\;\;\;\;\label{w40}\\
&&t\left({\bf d}^m_{MED}\right)=m-1\;,\;\;\;
\min\Delta_{{\cal H}}\left({\bf d}^m_{MED}\right)=d_m-d_{m-1}\;,\;\;\;\max
\Delta_{{\cal H}}\left({\bf d}^m_{MED}\right)=d_{m-1}-m\;,\nonumber\\
&&\#\Delta_{{\cal H}}\left({\bf d}^m_{MED}\right)=\frac{2}{m}\sum_{k=2}^md_k-
d_m\;,\;\;\;\#\Delta_{{\cal G}}\left({\bf d}^m_{MED}\right)=d_m-\frac{1}{m}
\sum_{k=2}^md_k-\frac{m-1}{2}\;.\nonumber
\end{eqnarray}
We need one more basic entity for ${\sf S}\left({\bf d}^m\right)$ which plays a 
key role and facilitates further discussion: {\em the Ap\'ery set} ${\mathbb A}
{\mathbb P}\left({\bf d}^m;d_1\right)$ of semigroup ${\sf S}\left({\bf d}^m
\right)$ with respect to generator $d_1$ is defined as follows,
\begin{eqnarray}
{\mathbb A}{\mathbb P}\left({\bf d}^m;d_1\right):=\left\{s\in {\sf S}\left(
{\bf d}^m\right)\quad |\quad s-d_1\not\in {\sf S}\left({\bf d}^m\right)\right\}
\;,\;\;\;\;\#{\mathbb A}{\mathbb P}\left({\bf d}^m;d_1\right)=d_1\;.
\label{ar14}
\end{eqnarray}
The generating function $AP_1\left({\bf d}^m;z\right)$ for the Ap\'ery set 
${\mathbb A}{\mathbb P}\left({\bf d}^m;d_1\right)$ was given in \cite{feai07}, 
Formula (4.4),
\begin{eqnarray}
AP_1\left({\bf d}^m;z\right)=\sum_{s\in {\mathbb A}{\mathbb P}\left({\bf d}^m;
d_1\right)}z^s=(1-z^{d_1})\cdot H\left({\bf d}^m;z\right)\;,\nonumber 
\end{eqnarray}
and is related to the numerator $Q\left({\bf d}^m;z\right)$ as follows, 
\begin{eqnarray}
Q\left({\bf d}^m;z\right)={\prod_{j=2}^m(1-z^{d_j})}\cdot AP_1\left({\bf d}^m;
z\right)\;.\label{w41}
\end{eqnarray}
As for the MED--semigroups, by \cite{heku71} we get
\begin{eqnarray}
{\mathbb A}{\mathbb P}\left({\bf d}^m_{MED};m\right)=\{0,d_2,\ldots ,d_m\}\;,\;
\;\;\;\;AP_1\left({\bf d}^m_{MED};z\right)=1+\sum_{k=2}^mz^{d_k}\;.\label{w42}
\end{eqnarray}
Now we arrive at the explicit expression for $Q\left({\bf d}^m_{MED};z\right)$
which is to our knowledge not discussed in literature. By insertion (\ref{w42}) 
into (\ref{w41}) we obtain,
\begin{eqnarray}
Q\left({\bf d}^m_{MED};z\right)=1+\sum_{k=1}^{m-2}(-1)^k\left[\;I_{m,k}(z)+
J_{m,k}(z)\;\right]+(-1)^{m-1}I_{m,m-1}(z)\;,\;\;\;\mbox{where}\;\;\;\;\;
\label{w43}\\
I_{m,k}(z)=\sum_{j_1>\ldots >j_{k-1}\geq 2}^mz^{2d_{j_1}+\overbrace{d_{j_2}+
\ldots+d_{j_{k}}}^{{\rm k-1}\;terms}},\;\;J_{m,k}(z)=k\sum_{j_1>\ldots >j_k
\geq 2}^mz^{\overbrace{d_{j_1}+d_{j_2}+\ldots+d_{j_{k+1}}}^{{\rm k+1}\;terms}}.
\;\label{w45}
\end{eqnarray}
The number of contributing monomials into $I_{m,k}(z)$ and $J_{m,k}(z)$ read
\begin{eqnarray}
\#I_{m,k}(z)=(m-1){m-2\choose k-1}\;,\;\;\;\;\;\#J_{m,k}(z)=k{m-1\choose k+1}\;.
\nonumber
\end{eqnarray}
Below we give the polynomials $I_{m,k}(z)$ and $J_{m,k}(z)$ for small and large 
indices $k$,
\begin{eqnarray}
&&I_{m,1}(z)=\sum_{j_1\geq 2}^mz^{2d_{j_1}},\;\;I_{m,2}(z)=\sum_{j_1>j_2\geq 2}
^mz^{2d_{j_1}+d_{j_2}},\;\;I_{m,3}(z)=\sum_{j_1>j_2>j_3\geq 2}^mz^{2d_{j_1}+
d_{j_2}+d_{j_3}},\;\ldots ,\nonumber\\
&&I_{m,m-2}(z)=\sum_{j_1>j_2>\ldots >j_{m-2}\geq 2}^mz^{2d_{j_1}+d_{j_2}+
\ldots+ d_{j_{m-2}}},\;\;\;\;I_{m,m-1}(z)=
z^{\Sigma_m}\sum_{j\geq 2}^mz^{d_j-m}\;,\label{w46}\\
&&J_{m,1}(z)=\sum_{j_1>j_2\geq 2}^mz^{d_{j_1}+d_{j_2}},\;J_{m,2}(z)=2
\sum_{j_1>j_2>j_3\geq 2}^mz^{d_{j_1}+d_{j_2}+d_{j_3}},\ldots ,\;
J_{m,m-2}(z)=(m-2)z^{\Sigma_m-m}\nonumber
\end{eqnarray}
In the presentation (\ref{w43}) it is easy to recognize the partial polynomials 
$Q_k\left({\bf d}^m_{MED};z\right)$ which are contributing to numerator 
$Q\left({\bf d}^m_{MED};z\right)$ in accordance with (\ref{bet05}),
\begin{eqnarray}
Q_k\left({\bf d}^m_{MED};z\right)=I_{m,k}(z)+J_{m,k}(z)\;,\;\;\;1\leq k\leq m-2
\;,\;\;\;Q_{m-1}\left({\bf d}^m_{MED};z\right)=I_{m,m-1}(z)\;.\label{w47}
\end{eqnarray}
The number of contributing terms $z^{C_{j,i}}$ into $Q_k\left({\bf d}^m_{MED};
z\right)$ coincides with $\beta_k\left({\bf d}^m_{MED}\right)$,
\begin{eqnarray}
\#Q_k\left({\bf d}^m_{MED};z\right)=\#I_{m,k}(z)+\#J_{m,k}(z)=(m-1){m-2\choose 
k-1}+k{m-1\choose k+1}=k{m\choose k+1}\;,\nonumber
\end{eqnarray}
in accordance with (\ref{w40}). Prove the main Theorem of this section.
\begin{theorem}\label{the7}Let a numerical MED--semigroup ${\sf S}\left({\bf d}
^m_{MED}\right)$ be given. Then it is almost symmetric iff for every element $d_
j$ of generating set ${\bf d}^m_{MED}$ there exists its counterpartner $d_{m-j+
1}$ such that
\begin{eqnarray}
d_j+d_{m-j+1}=m+d_m\;,\;\;\;\;\;1\leq j\leq m\;.\label{w48}
\end{eqnarray}
\end{theorem}
{\sf Proof} $\;\;\;$In accordance with (\ref{tr22}) a numerical MED--semigroup 
${\sf S}\left({\bf d}^m_{MED}\right)$ is almost symmetric iff 
\begin{eqnarray}
\Delta_{{\cal H}}\left({\bf d}^m_{MED}\right)=\left[{\mathbb B}_{m-1}\left(
{\bf d}^m_{MED}\right)\oplus\left\{-\Sigma_m\right\}\right]\setminus
\left\{F\left({\bf d}^m_{MED}\right)\right\}\;.\label{w48a}
\end{eqnarray}
According to (\ref{w46}) and (\ref{w47}) a set ${\mathbb B}_{m-1}\left({\bf d}^
m_{MED}\right)\oplus\left\{-\Sigma_m\right\}$ is composed of degrees of 
monomials $z^{d_j-m}$ entering the polynomial $I_{m,m-1}(z)$. In other words, 
by (\ref{w48a}) we have
\begin{eqnarray}
\Delta_{{\cal H}}\left({\bf d}^m_{MED}\right)=\left\{h_2,h_3,\ldots ,h_{m-1}
\right\}\;,\;\;\;h_j=d_j-m\;.\label{w49}
\end{eqnarray}
However, by definition (\ref{tr9}) of the set $\Delta_{{\cal H}}\left({\bf d}^m
\right)$ for every element $h_j\in\Delta_{{\cal H}}\left({\bf d}^m\right)$ there
exists its counterpartner $h_j^{\star}\in\Delta_{{\cal H}}\left({\bf d}^m
\right)$ such that $h_j+h_j^{\star}=F\left({\bf d}^m\right)$. Substituting the 
expression (\ref{w49}) for the gaps and the expression (\ref{w40}) for the 
Frobenius number into the last equality we come to the necessary and efficient 
conditions in the case of the almost symmetric MED-semigroup,
\begin{eqnarray}
d_j+d_j^{\star}=d_m+m\;,\;\;\;2\leq j\leq m-1\;.\label{w50} 
\end{eqnarray}
Since the tuple ${\bf d}^m_{MED}$ is arranged in accending order, $m<d_2<d_3<
\ldots<d_{m-1}<d_m$, then by (\ref{w50}) a set of counterpartners $d_j^{\star}$ 
has to be arranged in descending order, $d_2^{\star}>d_3^{\star}>\ldots >d_{m-1}
^{\star}$. Combining both sequences with opposite growth we can verify that 
(\ref{w50}) could be satisfied for every generator $d_j$ iff $d_j^{\star}=d_{m-
j+1}$. A proof can be given combining a way of contradiction with induction on 
index $j$ in (\ref{w48}). 

Since a case $j=1$ is trivial, we start with $j=2$. According to (\ref{w40}), 
(\ref{w48a}) and (\ref{w49}) we have
\begin{eqnarray}
d_2-m=\min\Delta_{{\cal H}}\left({\bf d}^m_{MED}\right)=d_m-d_{m-1}\;,\nonumber
\end{eqnarray}
that satisfy (\ref{w48}). Let equality (\ref{w48}) holds for all $1\leq j\leq 
q$. Prove, by way of contradiction, that it holds also for $j=q+1$. Indeed, let 
for elements $d_{q+1}\in{\bf d}^m_{MED}$ and $d_{m-q}\in{\bf d}^m_{MED}$ there 
exist counterpartners $d_{q+1}^{\star}\in{\bf d}^m_{MED}$ and $d_{m-q}^{\star}
\in{\bf d}^m_{MED}$, respectively, such that $d_{m-q}>d_{q+1}^{\star}$. In 
accordance with (\ref{w50}) write two equalities
\begin{eqnarray}
d_q+d_{m-g+1}=d_{m-q}+d_{m-q}^{\star}\;\;\;\;\;\mbox{and}\;\;\;\;\;
d_{q+1}+d_{q+1}^{\star}=d_{m-q}+d_{m-q}^{\star}\;,\nonumber
\end{eqnarray}
which give rise to following inequalities,
\begin{eqnarray}
d_{m-q}^{\star}-d_q=d_{m-g+1}-d_{m-q}>0\;,\;\;\;\;\;
d_{q+1}-d_{m-q}^{\star}=d_{m-q}-d_{q+1}^{\star}>0\;.\label{w49b}
\end{eqnarray}
Thus, by (\ref{w49b}) we arrive at inequality $d_q<d_{m-q}^{\star}<d_{q+1}$.  
However, the last inequality has not solutions because it presumes existence of 
generator $d_{m-q}^{\star}\in {\bf d}^m$ between $d_q$ and $d_{q+1}$ that 
contradicts the arrangement of the tuple $d_1<\ldots<d_q<d_{g+1}<\ldots<d_m$. 
This finishes proof of Theorem.$\;\;\;\;\;\;\Box$

Proposition \ref{pro1} comes as Corollary of Theorem \ref{the7}. Indeed, a tuple
of arithmetic sequence
$$
t+1,\;t+1+\frac{g}{t},\;t+1+2\frac{g}{t},\;\ldots ,\;t+1+g\;,
$$
with generic term $d_j=t+1+(j-1)\cdot g/t$ satisfies (\ref{w48}) : $d_j+d_{m-g+
1}=g+2(t+1)$.

Consider another Corollary which follows by Theorem \ref{the7}.
\begin{corollary}\label{cor3}Let an almost symmetric MED--semigroup ${\sf S}
\left({\bf d}^{2n+1}_{MED}\right)$ be given. Then an element $d_{2n+1}$ is 
an odd integer and a sum $\sum_{i=1}^{2n+1}d_i$ is divisible by $2n+1$.
\end{corollary}
{\sf Proof} $\;\;\;$Since $t\left({\bf d}^{2n+1}_{MED}\right)=2n$, then there 
exists an index $j=n+1$ such that $d_j=d_{2n+1-j+1}$ and by consequence of 
(\ref{w48}) the following equality holds, $2d_{n+1}=2n+1+d_{2n+1}$. Hence, it 
follows the 1st part of Corollary: $d_{2n+1}$ is an odd integer. The 2nd part 
follows if we denote, in accordance with the 1st part, $d_{2n+1}=2w+1$ and 
calculate,
\begin{eqnarray}
\sum_{i=1}^{2n+1}d_i=(d_1+d_{2n+1})+(d_2+d_{2n})+\ldots +(d_{n}+d_{n+2})+
d_{n+1}\;.\nonumber
\end{eqnarray}
According to Theorem \ref{the7} and the 1st part of this Corollary we obtain
\begin{eqnarray}
\sum_{i=1}^{2n+1}d_i=n(2n+1+2w+1)+n+1+w=(2n+1)(n+1+w)\;.\;\;\;\;\Box\label{w90}
\end{eqnarray}

Explicit formulas for the type $t\left({\bf d}^m_{MED}\right)$ and the Betti 
numbers $\beta_k\left({\bf d}^m_{MED}\right)$ give another opportunity to 
specify Theorems \ref{the5} and \ref{the6} in the case of almost symmetric 
MED--semigroups. 
Calculate a sum of the Betti numbers $\beta_k\left({\bf d}^m_{MED}\right)$ and 
check that it satisfies inequality (\ref{bet22}),
\begin{eqnarray}
\sum_{k=0}^{m-1}\beta_k\left({\bf d}^m_{MED}\right)=(m-2)\cdot 2^{m-1}+2\;.
\label{w51}
\end{eqnarray}
\begin{theorem}\label{the8}Let an almost symmetric MED--semigroup ${\sf S}\left(
{\bf d}^m_{MED}\right)$ be given. Then 
\begin{eqnarray}
\rho\left({\bf d}^m_{MED}\right)=\ell\left({\bf d}^m_{MED}\right)\;.\label{w52}
\end{eqnarray}
\end{theorem}
{\sf Proof} $\;\;\;$Keeping in mind $\rho\left({\bf d}^m\right)-\ell\left({\bf 
d}^m\right)=\delta\left({\bf d}^m\right)$, we prove Theorem for even and odd 
edim separately. First, consider an almost symmetric MED--semigroup ${\sf S}
\left({\bf d}^{2n}_{MED}\right)$ and calculate the sums of the Betti numbers 
$\beta_k\left({\bf d}^{2n}_{MED}\right)$ of even and odd indices separately.

\noindent
Keeping in mind $\beta_{2n-1}\left({\bf d}^{2n}_{MED}\right)=2n-1$ and making 
sum of (\ref{w51}) with (\ref{bet2}) we obtain,
\begin{eqnarray}
&&\beta_1\left({\bf d}^{2n}_{MED}\right)+\beta_3\left({\bf d}^{2n}_{MED}\right)+
\ldots +\beta_{2n-3}\left({\bf d}^{2n}_{MED}\right)=(n-1)\cdot (2^{2n-1}-2)\;,
\label{w53}\\
&&\beta_2\left({\bf d}^{2n}_{MED}\right)+\beta_4\left({\bf d}^{2n}_{MED}\right)
+\ldots +\beta_{2n-2}\left({\bf d}^{2n}_{MED}\right)=(n-1)\cdot 2^{2n-1}\;.
\nonumber
\end{eqnarray}
Applying now Theorem \ref{the5}, combine the 1st or the 2nd pairs of equalities 
in (\ref{w53}) and (\ref{jo4a}) and get
\begin{eqnarray}
(n-1)2^{2n-1}=2(n-1)4^{n-1}+\frac1{2}\;\delta\left({\bf d}^{2n}_{MED}\right)
\;\;\;\;\longrightarrow\;\;\;\;\delta\left({\bf d}^{2n}_{MED}\right)=0\;.
\label{w54}
\end{eqnarray}
Next, consider an almost symmetric MED--semigroup ${\sf S}\left({\bf d}^{2n+1}_
{MED}\right)$ and make similar calculations,
\begin{eqnarray}
&&\beta_1\left({\bf d}^{2n+1}_{MED}\right)+\beta_3\left({\bf d}^{2n+1}_{MED}
\right)+\ldots +\beta_{2n-1}\left({\bf d}^{2n+1}_{MED}\right)=(2n-1)\cdot 
2^{2n-1}+1\;,\label{w55}\\
&&\beta_2\left({\bf d}^{2n+1}_{MED}\right)+\beta_4\left({\bf d}^{2n+1}_{MED}
\right)+\ldots +\beta_{2n-2}\left({\bf d}^{2n+1}_{MED}\right)=(2n-1)\cdot 2^{
2n-1}-2n\;.\nonumber
\end{eqnarray}
Applying Theorem \ref{the6}, combine the 1st or the 2nd pairs of equalities
in (\ref{w55}) and (\ref{jo5c}) and get,
\begin{eqnarray}
(2n-1)\cdot 2^{2n-1}-2n=(2n-1)\left(2^{2n-1}-1\right)-1+\frac1{2}\;\delta\left(
{\bf d}^{2n+1}_{MED}\right)\;\;\;\;\longrightarrow\;\;\;\;\delta\left({\bf d}^
{2n+1}_{MED}\right)=0\;.\label{w56} 
\end{eqnarray}
Thus, combining (\ref{w54}) and (\ref{w56}) we arrive at (\ref{w52}).
$\;\;\;\;\;\;\Box$

We finish this section with very specific almost symmetric MED--semigroup ${\sf 
S}\left({\bf d}^m_{MED}\right)$, related to Proposition \ref{pro1}, when both 
cardinalities $\rho\left({\bf d}^m_{MED}\right)$ and $\ell\left({\bf d}^m_{MED}
\right)$ are vanishing. It enhances an equality (\ref{w52}) in Theorem 
\ref{the8} in this specific case. First, we start with auxiliary Lemma.
\begin{lemma}\label{lem8}Let an almost symmetric semigroup ${\sf S}\left({\bf
d}^{2n+1}\right)$ be given and let the sets ${\mathbb B}_i\left({\bf d}^{2n+1}
\right)$ and ${\overline{\mathbb B}}_i\left({\bf d}^{2n+1}\right)$ be defined in
(\ref{ar10a}) and (\ref{z35}). Define four union sets,
\begin{eqnarray}
{\mathbb B}_o\left({\bf d}^{2n+1}\right)=\cup_{i=1}^n{\mathbb B}_{2i-1}\left(
{\bf d}^{2n+1}\right)\;,\;\;\;\;{\mathbb B}_e\left({\bf d}^{2n+1}\right)=\cup_{
i=1}^{n-1}{\mathbb B}_{2i}\left({\bf d}^{2n+1}\right)\;,\label{z60}\\
{\overline{\mathbb B}}_o\left({\bf d}^{2n+1}\right)=\cup_{i=1}^{n}{\overline
{\mathbb B}}_{2i-1}\left({\bf d}^{2n+1}\right)\;,\;\;\;\;{\overline{\mathbb B}}
_e\left({\bf d}^{2n+1}\right)=\cup_{i=1}^{n-1}{\overline{\mathbb B}}_{2i}\left(
{\bf d}^{2n+1}\right)\;.\nonumber
\end{eqnarray}
If a set ${\mathbb Y}=\left[{\mathbb B}_o\left({\bf d}^{2n+1}\right)\cup{
\overline{\mathbb B}}_e\left({\bf d}^{2n+1}\right)\right]\cap\left[{\mathbb B}_
e\left({\bf d}^{2n+1}\right)\cup{\overline{\mathbb B}}_o\left({\bf d}^{2n+1}
\right)\right]$ is empty then $\rho\left({\bf d}^{2n+1}\right)=0$.
\end{lemma}
{\sf Proof} $\;\;\;$Consider a multiset ${\mathfrak X}_{12}\left({\bf d}^{2n+1}
\right)$ defined in (\ref{gp5}) and write its standard representation
\begin{eqnarray}
{\mathfrak X}_{12}\left({\bf d}^{2n+1}\right)=\la{\mathbb X}_{12}\left({\bf d}^
{2n+1}\right),\sigma_{{\mathbb X}_{12}\left({\bf d}^{2n+1}\right)}\ra\;,
\label{z61}
\end{eqnarray}
where in view of definitions (\ref{bi3}) and (\ref{bi3a}) of multiset operations
$\bigvee$ and $\bigwedge$ the underlying set ${\mathbb X}_{12}\left({\bf d}^{2n
+1}\right)$ is given by
\begin{eqnarray}
{\mathbb X}_{12}\left({\bf d}^{2n+1}\right)=\left[{\mathbb B}_o\left({\bf d}^{
2n+1}\right)\cup{\overline{\mathbb B}}_e\left({\bf d}^{2n+1}\right)\right]\cap
\left[{\mathbb B}_e\left({\bf d}^{2n+1}\right)\cup{\overline{\mathbb B}}_o\left(
{\bf d}^{2n+1}\right)\right]\;.\label{z62}
\end{eqnarray}
By comparison of the sets ${\mathbb Y}$ and ${\mathbb X}_{12}\left({\bf d}^{2n+
 1}\right)$ we conclude that they concide, ${\mathbb Y}\equiv{\mathbb X}_{12}
\left({\bf d}^{2n+1}\right)$. However, by definition (\ref{z7}) of empty 
multiset a set equality ${\mathbb X}_{12}\left({\bf d}^{2n+1}\right)=\emptyset$ 
implies a multiset equality ${\mathfrak X}_{12}\left({\bf d}^{2n+1}\right)=
{\mathfrak \not O}$. Thus, $\rho\left({\bf d}^{2n+1}\right)=0$ and Lemma is 
proven.$\;\;\;\;\;\;\Box$

Denote by ${\bf c}$ the MED-tuple of edim=$2n+1$ such that its generating set 
is arranged as an arithmetic sequence and, according to Theorem \ref{the7},  
generates an almost symmetric semigroup,
\begin{eqnarray}
{\bf c}=\{2n+1,\;2n+1+2a,\;\ldots,\;2n+1+4na\}\;,\;\;\;a\in{\mathbb N}\;,\;\;
\gcd(2n+1,a)=1\;.\label{w58}    
\end{eqnarray}
\begin{corollary}\label{cor4}Let an almost symmetric MED--semigroup ${\sf S}
\left({\bf c}\right)$ defined in (\ref{w58}) be given. Then
\begin{eqnarray}
\rho\left({\bf c}\right)=\ell\left({\bf c}\right)=0\;.\label{w57}
\end{eqnarray}
\end{corollary}
{\sf Proof} $\;\;\;$We analyze a set ${\mathbb X}_{12}\left({\bf d}^{2n+1}
\right)$ defined in (\ref{z62}) for the case ${\bf d}^{2n+1}={\bf c}$ given in 
(\ref{w58}). Observe that elements of corresponding sets ${\mathbb B}_o\left(
{\bf c}\right)$, ${\overline{\mathbb B}}_o\left({\bf c}\right)$ and ${\mathbb B}
_e\left({\bf c}\right)$, ${\overline{\mathbb B}}_e\left({\bf c}\right)$ defined 
in (\ref{z60}) through the partial sets ${\mathbb B}_k\left({\bf c}\right)$ and 
${\overline{\mathbb B}}_k\left({\bf c}\right)$ are coming as degrees $\xi_{q,k}$
and $\theta_{q,k}$ of monomials $z^{\xi_{q,k}}\in I_{2n+1,k}(z)$ and $z^{\theta
_{q,k}}\in J_{2n+1,k}(z)$, respectively, and as their conjugates ${\overline 
\xi}_{q,k}$ and ${\overline \theta}_{q,k}$,
\begin{eqnarray}
\left\{\begin{array}{r}{\overline \xi}_{q,k}\\{\overline \theta}_{q,k}\end
{array}\right\}=-\left\{\begin{array}{r}\xi_{q,k}\\\theta_{q,k}\end{array}
\right\}+F\left({\bf c}\right)+\sum_{j=1}^{2n+1}d_j\;,\;\;\;\;\;z^{\xi_{q,k}}
\in I_{2n+1,k}(z)\;,\;\;\;z^{\theta_{q,k}}\in J_{2n+1,k}(z)\;.\;\;\;\;
\label{w72}
\end{eqnarray}
Indeed, by (\ref{w47}) every partial sets ${\mathbb B}_k\left({\bf c}\right)$ 
and ${\overline{\mathbb B}}_k\left({\bf c}\right)$ can be decomposed in other 
two sets,
\begin{eqnarray}
{\mathbb B}_k\left({\bf c}\right)={\mathbb B}_k^I\left({\bf c}\right)\cup
{\mathbb B}_k^J\left({\bf c}\right)\;,\;\;\;\;\;\;{\overline{\mathbb B}}_k\left(
{\bf c}\right)={\overline{\mathbb B}}_k^I\left({\bf c}\right)\cup{\overline
{\mathbb B}}_k^J\left({\bf c}\right)\;,\;\;\;\;\;\mbox{where}\label{w73a}
\end{eqnarray}
\begin{eqnarray}
{\mathbb B}_k^I\left({\bf c}\right)=\cup_{q=1}^{\beta_k}\{\xi_{q,k}\}\;,\;\;\;
{\mathbb B}_k^J\left({\bf c}\right)=\cup_{q=1}^{\beta_k}\{\theta_{q,k}\}\;,\;\;
\;{\overline{\mathbb B}}_k^I\left({\bf c}\right)=\cup_{q=1}^{\beta_k}\{
{\overline\xi}_{q,k}\}\;,\;\;\;{\overline{\mathbb B}}_k^J\left({\bf c}\right)=
\cup_{q=1}^{\beta_k}\{{\overline\theta}_{q,k}\}\;.\;\;\label{w73}
\end{eqnarray}
Consider parity properties of these elements. First, note that according to 
Corollary \ref{cor3} and (\ref{w90}) the following sum always takes odd values,
\begin{eqnarray}
F\left({\bf c}\right)+\sum_{j-1}^{2n+1}d_j=4an+(2n+1)(2n(a+1)+1)\;.\nonumber
\end{eqnarray}
The last equality together with (\ref{w72}) results in important conclusion: 
\begin{eqnarray}
\mbox{Elements}\;\;\;\xi_{q,k}\;\;\mbox{and}\;\;{\overline \xi}_{q,k}\;\;\;\;
\mbox{are of opposite parities as well as elements}\;\;\;\theta_{q,k}\;\;\mbox{
and}\;\;{\overline \theta}_{q,k}\;.\label{w74}
\end{eqnarray}
Consider the elements $\xi_{q,k}$ and $\theta_{q,k}$ in more details. By 
(\ref{w45}) write them as follows,
\begin{eqnarray}
\xi_{q,k}=2d_{j_1}+\overbrace{d_{j_2}+\ldots+d_{j_{k}}}^{{\rm k-1}\;terms}\;,
\;\;\;\;\;\;\theta_{q,k}=\overbrace{d_{j_1}+d_{j_2}+\ldots+d_{j_{k+1}}}^{{\rm 
k+1}\;terms}\;,\label{w75}
\end{eqnarray}
and recall that a generic term of the sequence (\ref{w58}) reads $d_j=2n+1+
2(j-1)a$. Combining it with (\ref{w74}) and (\ref{w75}) we conclude 
\begin{eqnarray}
&&2\mid \xi_{q,2k+1}\;,\;\;2\mid \theta_{q,2k+1}\;,\;\;\;\;\mbox{and}\;
\;\;\;2\nmid \xi_{q,2k}\;,\;\;2\nmid \theta_{q,2k}\;,\label{w76}\\
&&2\nmid {\overline \xi}_{q,2k+1}\;,\;\;2\nmid {\overline \theta}_{q,2k+1}\;,
\;\;\;\;\mbox{and}\;\;\;\;2\mid {\overline \xi}_{q,2k}\;,\;\;2\mid {\overline 
\theta}_{q,2k}\;.\nonumber
\end{eqnarray}
Thus, by (\ref{w73}) the sets ${\mathbb B}_{2k+1}^I\left({\bf c}\right)$, 
${\mathbb B}_{2k+1}^J\left({\bf c}\right)$ and ${\overline{\mathbb B}}_{2k}^I
\left({\bf c}\right)$, ${\overline{\mathbb B}}_{2k}^J\left({\bf c}\right)$ 
comprise the elements divisible by 2, while the sets ${\mathbb B}_{2k}^I\left(
{\bf c}\right)$, ${\mathbb B}_{2k}^J\left({\bf c}\right)$ and ${\overline
{\mathbb B}}_{2k+1}^I\left({\bf c}\right)$, ${\overline{\mathbb B}}_{2k+1}^J
\left({\bf c}\right)$ comprise the elements nondivisible by 2.

Next, based on the last conclusion and equalities (\ref{z60}) and (\ref{w73a}) 
we arrive at parity properties:
\begin{eqnarray}
&&\mbox{Sets}\;\;{\mathbb B}_o\left({\bf c}\right)\;\;\mbox{and}\;\;{\overline
{\mathbb B}}_e\left({\bf c}\right)\;\;\mbox{comprise only the elements divisible
by 2}\;,\;\;\;\;\;\;\label{w77}\\
&&\mbox{Sets}\;\;{\overline{\mathbb B}}_o\left({\bf c}\right)\;\;\mbox{and}\;\;
{\mathbb B}_e\left({\bf c}\right)\;\;\mbox{comprise only the elements 
nondivisible by 2}\;.\;\;\;\;\;\;\nonumber
\end{eqnarray}
Finally, according to (\ref{z62}) the set ${\mathbb X}_{12}\left({\bf c}\right)$
is empty since, by (\ref{w77}), two pairs of sets, ${\mathbb B}_o\left({\bf c}
\right)\cup {\overline{\mathbb B}}_e\left({\bf c}\right)$ and ${\mathbb B}_e
\left({\bf c}\right)\cup {\overline{\mathbb B}}_o\left({\bf c}\right)$, comprise
elements of distinct parities. Therefore, by Lemma \ref{lem8} this implies $\rho
\left({\bf c}\right)=0$. However, by Theorem \ref{the8} the last equality leads 
immediately to another equality, $\ell\left({\bf c}\right)=0$, that finishes 
our proof.$\;\;\;\;\;\;\Box$
\section*{Acknowledgement}
The useful discussions with A. Juhasz are highly appreciated. 

\end{document}